\documentclass[
  a4paper, 
  reqno, 
  oneside, 
  11pt
]{amsart}

\usepackage[utf8]{inputenc}

\usepackage[dvipsnames]{xcolor} 
\colorlet{cite}{LimeGreen!50!Green}
\usepackage{tikz}

\usetikzlibrary{arrows,positioning, calc, 3d, perspective, arrows.meta}
\tikzset{ 
  baseline=-2.3pt,
  text height=1.5ex, text depth=0.25ex,
  >=stealth,
  node distance=2cm,
  mid/.style={fill=white,inner sep=2.5pt},
}

\usepackage{lmodern}
\usepackage{tikz-cd}

\usetikzlibrary{
matrix,
arrows,
calc,
intersections,
decorations.markings,
decorations.pathmorphing,
shapes,
positioning,
}

\usepackage{amsthm, amssymb, amsfonts}
\usepackage{tabularray}
\usepackage{cancel}
\usepackage{graphicx,caption,subcaption}
\usepackage{braket}
\usepackage{microtype}

\usepackage[%
  bookmarks=true,     
  unicode=true,     
  pdftitle={},    %
  pdfauthor={Elizabeth Gasparim}{Bruno Suzuki}, %
  pdfkeywords={}, 
  colorlinks=true,    
  linkcolor=Blue,     
  citecolor=cite,   
  filecolor=magenta,    
  urlcolor=RoyalBlue      
]{hyperref}       
\usepackage{cleveref}

\newtheoremstyle{mydef}
  {}
  {}
  {}
  {}
  {\scshape}
  {. }
  { }
  {\thmname{#1}\thmnumber{ #2}\thmnote{ #3}}

\theoremstyle{plain}
\newtheorem{theorem}{Theorem}

\newtheorem*{theorem*}{Theorem}

\theoremstyle{mydef}
\newtheorem{definition}[theorem]{Definition}
\newtheorem*{conjecture*}{Conjecture}
\theoremstyle{remark}

\newtheorem{notation}[theorem]{Notation}
\newtheorem{example}[theorem]{Example}

\newtheorem*{proposition*}{Proposition}
\newtheorem*{lemma*}{Lemma}
\newtheorem*{corollary*}{Corollary}
\theoremstyle{definition}
\theoremstyle{remark}

\DeclareMathOperator{\Tot}{Tot}

\newcommand{\ce}{\mathrel{\mathop:}=}
\newcommand{\multitilde}[2]{\widetilde{#2}^{\raisebox{0.3ex}{$\scriptscriptstyle#1$}}}
\newcommand{\newscale}{0.5}

\usepackage[dvipsnames]{xcolor}
\colorlet{cite}{LimeGreen!50!Green}

\newcommand{\stA}{[RoyalBlue,thick]}
\newcommand{\stB}{[thick,dashed]}

\begin{document}

\author{Elizabeth Gasparim}
\address{E.G. UTFPR, Brazil}
\email{etgasparim@gmail.com}

\author{Carlos Simpson}
\address{C.S. CNRS Université Côte d'Azur, France}
\email{carlos.simpson@univ-cotedazur.fr}

\author{Bruno Suzuki}
\address{B.S. UNESPAR, Brazil}
\email{obrunosuzuki@gmail.com}

\author{Rashid Talha}
\address{R.T. ICTP, Italy}
\email{rtalha@ictp.it}

\title[Calabi--Yau singularities and quilts]{Calabi--Yau threefold singularities\\
and their universal quilts}

\begin{abstract}
We introduce the notion of universal quilt for toric CY3 singularities and construct examples when the universal quilts are themselves smooth Calabi--Yau threefolds. Such quilts are designed to contain all of the crepant resolutions of the given singularity. To obtain universal quilts we implement an algorithm that computes all triangulations for any given planar polygon. We also introduce the notion of quilt stack, to address the question of finding the most likely smooth configuration preceding a given singularity.
\end{abstract}

\maketitle

\tableofcontents

\section{Crepant Resolutions of Calabi--Yau Singularities}
We invite the reader to go first of all to subsection \ref{pictures} in the hope that 
seeing those pictures will serve  as a geometric motivation. 
Besides drawing artistic pictures, our motivation for studying crepant resolutions of Calabi--Yau singularities is dual. 
Firstly, crepant resolutions  are needed  when computing various partition functions  in String Theory
when the target space is singular, e.g.  \cite{Ne1,Ne2, NO, GL, GV1,GV2}.
Secondly, for questions in Cosmology there is significant  interest
in  predicting, out of all possible scenarios, which 
one is most likely to have preceded a given singularity; thus, estimating the past configuration of such a singularity.

In the string theory case, typically, results are known for  smooth targets, 
and given a singular target, authors will choose one of its resolutions.
In our opinion, such an approach will only bring to light 
features of the particular choice of resolution, instead of those of the original target. 
In \cite{GSTV, GKMR} a method  for computing partition functions for singular varieties was proposed, which does not
restrict to a fixed choice of resolution, instead it uses an average of all possible resolutions;
thus,  requiring knowledge of all possible resolutions. \\

Finding all crepant resolutions of a toric Calabi--Yau threefold leads us  to  a very concrete hard problem in combinatorics \cite{Eps}; namely, that of finding all primitive triangulations of a given polytope. It is well known that a toric Calabi--Yau threefold is necessarily noncompact. Furthermore, a toric variety $X$ is Calabi--Yau if and only if all of the vectors its defining fan end up on the same hyperplane. In particular, for the 3 dimensional case, this means that the heads $p_1, \dots, p_n$ of the vectors $v_i=\overrightarrow{0p_i}$ forming the fan of $X$ must all end on a plane (here the fan is being described by vectors starting at the origin). Consequently, the Calabi--Yau threefold $X$ itself is fully determined by the  planar polytope $P$ with vertices $v_1, \dots, v_n$. We will always work with an integer lattice (the lattice of characters of the torus) and with polytopes whose vertices belong to this integer lattice.

In general, $P$ determines a singular Calabi--Yau threefold. Suppose $P$ is a lattice polytope with area at least $1$. Then it means that $P$ contains at least $4$  lattice points, say $p_1, \dots, p_4$. Accordingly, the vectors $v_i=\overrightarrow{0p_i}$ are linearly dependent in $\mathbb C^3$, and this implies that the corresponding variety is cut out by an algebraic equation, determining a singular threefold contained in $\mathbb C^4$. Typically one thinks of the planar polytope as being positioned at the hyperplane $z=1$. As an example, for the square below, we  write the vectors of the fan of $X$ as $v_1=(0,0,1), v_2=(1,0,1), v_3=(1,1,1), v_4=(0,1,1)$. The inward pointing normals $m_i$ of the facets of the fan determine the equations of the variety. In this case, the normals satisfy an equation of the form $m_1+m_2=m_3+m_4$, therefore the corresponding torus characters $x_i = \chi^{m_i}$ satisfy the equation $x_1x_2-x_3x_4=0$; which is the equation cutting out the singular threefold in $\mathbb C^4$. It is true in general that  toric varieties are cut out by such binomial equations. For more details on how the inward normals are calculated and how the equation of the Calabi--Yau singularity is found from its planar polygon see Example \ref{strip24}. For the underlying theory of toric varieties see the foundational lecture notes of David Cox available in his homepage.\footnote{Here is a video of dual fans in 3d \url{https://www.youtube.com/watch?v=XJ-WmmaJPAY}.}

\begin{center}

\begin{tikzpicture}[scale=0.75]
\draw (0,0) -- (1,0) -- (1,1) -- (0,1) -- cycle;
\node[xshift=-0.08cm, yshift=-0.08cm] at (0,0) {$p_1$};
\node[xshift=0.25cm, yshift=-0.08cm] at (1,0) {$p_2$};
\node[xshift=0.25cm, yshift=0.25cm] at (1,1) {$p_3$};
\node[xshift=-0.08cm, yshift=0.25cm] at (0,1) {$p_4$};
\draw[green] (0,0) -- (1,1);
\end{tikzpicture}
\end{center}

Getting back to the example of the square,  
the convex hull of $p_1, \dots, p_4$ can then be subdivided by  an additional edge, say $p_1p_3$, depicted in green. 
In the geometry of the threefold, the added edge of the polygon represents a $\mathbb P^1$ obtained by doing a small blowup (a crepant resolution of $X$).
Now, the triangle $p_1p_2p_3$ determines a basis $v_1v_2v_3$ of $\mathbb C^3$,
so that the fan $\{v_1,v_2,v_3\}$ describes a smooth threefold, and similarly for $p_1p_3p_4$. 
More vertices on $P$ just imply a larger number of triangles, all of which should be of area $1/2$  for obtaining a smooth variety. 
In conclusion, each crepant resolution of $X$ corresponds to a primitive triangulation of $P$, that is, formed by triangles of area $1/2$. 
For details about toric resolution of singularities see the very clear text of D. Cox  \cite[\S 5]{Cox}. \\
 
\begin{figure}
\centering
\definecolor{beige}{rgb}{0.96, 0.96, 0.86}
\definecolor{carnationpink}{rgb}{1.0, 0.65, 0.79}
\definecolor{blush}{rgb}{0.87, 0.36, 0.51}
\newcommand{\newnewcolor}{blush}
\newcommand{\fanscale}{1.2}
\begin{tikzpicture}[3d view={155}{30}, scale=\fanscale] 

\path	coordinate (O) at (0,0,0)
		coordinate (A) at (0,0,2)
		coordinate (B) at (4,0,2)
		coordinate (C) at (2,1,2)
		coordinate (D) at (0,1,2);

\draw[] (O) -- +(3.5,0,0);
\draw[] (O) -- +(0,1.5,0);
\draw[] (O) -- +(0,0,2.5);

\draw[fill = beige] (O) -- (A) -- (B) -- cycle;
\draw[fill = lightgray] (O) -- (B) -- (C) -- cycle;

\draw[fill = beige] (O) -- (D) -- (A) -- cycle;

\draw[thick, black] (A) -- (B) -- (C) -- (D) -- cycle;

\draw[fill=black] (0,0,2) circle (0.03);
\draw[thick, -{Latex[length=2mm]}, \newnewcolor] (O) -- (A);

\draw[fill = lightgray] (O) -- (C) -- (D) -- cycle;

\draw[fill=black] (0,1,2) circle (0.03);
\draw[fill=black] (1,1,2) circle (0.03);
\draw[fill=black] (2,1,2) circle (0.03);
\draw[fill=black] (1,0,2) circle (0.03);
\draw[fill=black] (2,0,2) circle (0.03);
\draw[fill=black] (3,0,2) circle (0.03);
\draw[fill=black] (4,0,2) circle (0.03);

\draw[thick, -{Latex[length=2mm]}, \newnewcolor] (O) -- (B);

\draw[thick, -{Latex[length=2mm]}, \newnewcolor] (O) -- (D);

\draw[thick, -{Latex[length=2mm]}, \newnewcolor] (O) -- (C);

\end{tikzpicture}   \quad
\begin{tikzpicture}[3d view={155}{30}, scale=\fanscale] 

\path	coordinate (O) at (0,0,0)
		coordinate (A) at (0,0,2)
		coordinate (B) at (4,0,2)
		coordinate (C) at (2,1,2)
		coordinate (D) at (0,1,2);

\draw[] (O) -- +(3.5,0,0);
\draw[] (O) -- +(0,1.5,0);
\draw[] (O) -- +(0,0,2.5);

\draw[fill = beige] (O) -- (A) -- (B) -- cycle;
\draw[fill = lightgray] (O) -- (B) -- (C) -- cycle;

\draw[fill = beige] (O) -- (D) -- (A) -- cycle;

\draw[thick, black] (A) -- (B) -- (C) -- (D) -- cycle;

\draw[fill=black] (0,0,2) circle (0.03);
\draw[thick, -{Latex[length=2mm]}, \newnewcolor] (O) -- (A);

\draw[fill = lightgray] (O) -- (C) -- (D) -- cycle;

\draw[fill=black] (0,1,2) circle (0.03);
\draw[fill=black] (1,1,2) circle (0.03);
\draw[fill=black] (2,1,2) circle (0.03);
\draw[fill=black] (1,0,2) circle (0.03);
\draw[fill=black] (2,0,2) circle (0.03);
\draw[fill=black] (3,0,2) circle (0.03);
\draw[fill=black] (4,0,2) circle (0.03);

\draw[thick, -{Latex[length=2mm]}, \newnewcolor] (O) -- (B);

\draw[thick, -{Latex[length=2mm]}, \newnewcolor] (O) -- (D);

\draw[thick, -{Latex[length=2mm]}, \newnewcolor] (O) -- (C);

\draw[thick, dashed] (1,1,2) -- (0,0,2);
\draw[thick, dashed] (1,1,2) -- (1,0,2);
\draw[thick, dashed] (2,1,2) -- (1,0,2);
\draw[thick, dashed] (2,1,2) -- (2,0,2);
\draw[thick, dashed] (2,1,2) -- (3,0,2);
\end{tikzpicture}   
\caption*{The fan of a toric CY3 singularity and one  of its resolutions}\label{CYfan}
\end{figure}
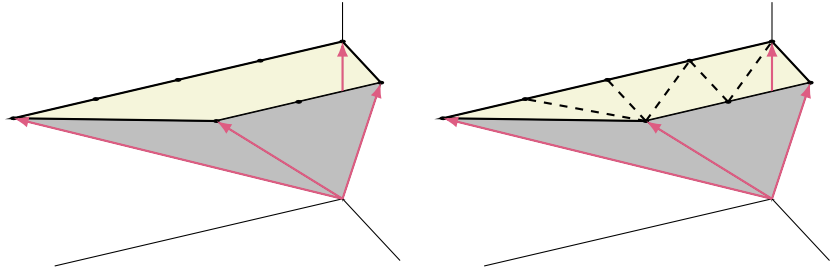

In this work, we discuss the geometric background underlying the set of  all crepant resolutions of a singular Calabi-Yau threefold. 
Focusing on a list of 16 examples of polytopes each containing 1 lattice point in their interior
(section \ref{favpolys}) we prove the following results:

\begin{theorem*}[\ref{tri}]
Let $X$ be a toric Calabi--Yau threefold singularity whose planar polytope 
is a convex triangle. Then it admits a smooth toric Calabi--Yau threefold as a universal quilt $\mathcal U(X)$. 
\end{theorem*}

\begin{theorem*}[\ref{par}]
Let $X$ be a toric Calabi--Yau threefold singularity whose planar polytope 
is a parallelogram. Then it admits a smooth toric Calabi--Yau threefold as a universal quilt $\mathcal U(X)$. 
\end{theorem*}

We observe that a fixed polygon may have more than one convex quilt, see Examples \ref{chess} and \ref{scarf}, 
while it is also possible for a polygon to admit various quilts, none of which are convex, see Example \ref{ali}.

\begin{theorem*}[\ref{tah}] (Tahzibi\footnote{We thank Ali Tahzibi for giving us a proof of this result.})
Let $X$ be a toric Calabi--Yau threefold singularity whose planar polytope 
is a  hexagon. Then the quilts of $X$ do not represent toric varieties. 
\end{theorem*}

For a chosen polygon, having drawn all of its  primitive triangulations (Def.~\ref{primi}), we can then consider 
the set of all such triangulations modulo $SL^{\pm}(2,\mathbb Z)$ transformations. We call such a 
quotient a quilt stack (Def.~\ref{qstack}).
Then 2 possibilities occur, and their studies require different approaches, they are:
\begin{itemize}
\item The triangulations are all pairwise distinct and the stack has height 1.  In such a case, we discuss 
chopping up corners of the polygon (these denote  toric blowups) and given that such blowups 
require energy to be accomplished,  we
declare the probability of a configuration as being the highest when it requires least blowing ups. 
Such is the case of all triangulations of polygon 13 studied in section \ref{toy}.
\item The triangulations are divided up into $SL^{\pm}(2,\mathbb Z)$ classes, at least one of which contain more 
than one element. When the largest equivalence class contains $n$ elements, then we obtain a quilt stack of 
height $n$. See Examples \ref{stack1} and \ref{stack2} for stacks of heights 8 and 6. 
\end{itemize}

Let us finish the introduction by mentioning some very interesting literature and social pages related to the work we present here.

Colleagues who were around when we printed out some figures, asked us
whether we would sew actual physical quilts made of cloth.  We are aware of \cite[\S 9]{IS},
where Irena Swanson discusses the mathematics and the making of semiregular tessellation quilts; 
we challenge quilters to produce in cloth our versal or universal quilts.

Sir Roger Penrose created in the 70's  the first examples of the famous figures known as Penrose Tiles \cite{RP1,RP2}.
Our universal quilts appear not to be examples of Penrose tiles, in the sense that each subtriangulation occurs a unique time in the quilt.

We emphasise that each universal quilt contains all possible primitive triangulations of a given polygon. Equivalently, 
a universal quilt contains all possible crepant resolutions of singularities of the initial CY3 singularity. 
While in this first paper on universal quilts we only carry out the construction  for a few examples, 
our algorithm applies to any convex polytope, therefore we are presenting  tools for making infinitely many universal quilts.

\section{Description of the Triangulations Algorithm}
The enumeration of triangulations of planar polygons is a hard problem in computational geometry, and its complexity depends strongly on the class of point sets under consideration \cite{Eps}.
Motivated by the need to enumerate the triangulations of convex lattice polygon as discussed in \S 1, we developed a specialised Python script, \texttt{triangulations.py}, for computing all primitive triangulations of a given convex lattice polygon. The code is hosted on GitHub\footnote{\url{https://github.com/rashidtalha/polytope-triangulations}} along with various additional scripts for visualisation.

Our approach is inspired by the method of Ray and Seidel for triangulations of planar point sets \cite{RaySeidel}. However, the structure of the present problem considerably simplifies the implementation. Since all polygons lie on a 2d lattice and we only need primitive triangulations, the number of edges in any triangulation is fixed. As a result, the problem can therefore be reduced to the enumeration of suitable non-intersecting subsets of candidate edges, which we perform using a depth-first search. This is closely related to the discussion on enumeration and counting of triangulations developed in, for example, \cite{AlvarezSeidel,Bespamyatnikh}.

The user provides an ordered list $S=(v_{1},\ldots,v_{k})$ of vertices corresponding to a convex lattice polygon $P=\overline{\mathrm{conv}(S)}$. The algorithm then computes all primitive triangulations of the lattice point set $P\cap\mathbb{Z}^{2}$. Each triangulation is returned as a list of its internal edges. The origin of coordinates can be chosen arbitrarily, but we require the vertices to be listed in counter-clockwise order, by convention.

In principle, $S$ can be arbitrarily large, although the runtime increases significantly with the number of lattice points in $P$, and in particular with the number of interior lattice points.

The algorithm consists of four main stages.

\subsection{Boundary and Internal Points}
Since $S$ only specifies the vertices of $P$, we start by computing the boundary points $\partial P \cap \mathbb{Z}^{2}$ by traversing the given vertices and interpolating intermediate integer coordinates. We also compute all strictly interior lattice $P^{\circ} \cap \mathbb{Z}^{2}$ points by scanning a bounding box and validating point containment via half-plane intersection tests using the usual cross product checks. Let $I$ and $B$ denote the number of boundary points and internal points, respectively.

\subsection{Candidate Edges}
Next, the algorithm generates a complete list of allowed internal edges (diagonals). A line segment connecting two lattice points is a possible candidate if and only if its interior contains no other lattice points. This is ensured by requiring that the GCD of the absolute coordinate differences of the endpoints is exactly $1$. Boundary edges are excluded from this candidate pool, as they always be included in every triangulation.

\subsection{Intersection Matrix}
In a valid triangulation, no two edges can intersect. So, the algorithm computes a Boolean matrix $M$ which shows whether any given pair of candidate edges intersect or not. This will play a central role in the final stage of the algorithm.

By Pick's theorem, the total area of $P$ is $A = I + B/2 - 1$. Therefore, each triangulation of $P$ will contain exactly $T = 2I + B - 2$ triangles of area $1/2$. Furthermore, by applying the Euler's formula $V - E + F = 2$ to the planar graph formed by the triangulation, with $V = I + B$ and $F = T + 1$ (the $T$ bounded faces plus the one unbounded face), we conclude that $E = 3I + 2B - 3$. There are exactly $B$ boundary edges, so the number of internal edges in each triangulation of $P$ is $E_{i} = 3I + B - 3$.

In particular, $E_{i}$ is completely determined by $P$, and this invariant significantly reduces the subsequent search space.

\subsection{Depth-First Search on Valid Edges}
Finally, the algorithm using the standard depth-first search (DFS) to find all subsets of candidate edges of size exactly $E_{i}$ that are pairwise non-intersecting (using matrix $M$ as a lookup table). Since DFS explores all such subsets of size $E_{i}$, it ensure that no triangulation is missed.

As an additional step, the algorithm can optionally compute a planar graph embedding for each valid configuration. It constructs an adjacency list, sorts the neighbours of each vertex according to their polar angle, and extracts the triangular faces by traversing the resulting cycles. This step is used only for subsequent visualisation and contructing the associated quilts.

\subsection{Reinforcement learning for the algorithm}
We note that one can apply reinforcement learning to the choice of next triangles in the search phase of the algorithm along the lines of \cite{semigr}. This is not implemented in the current version of the code. It should allow to treat triangulations of larger polytopes, for example the $3\times 3$ square with $4$ internal points has 
$46456$ triangulations and an $m\times n$ rectangle has very roughly on the order of $4^{mn}$ triangulations, see \cite{KZ}.

\section{Reading the Calabi-Yau Threefolds out of the Polygons}

It is well-known that in toric geometry the figure of a triangle represents $\mathbb P^2$ 
while  the figure of a square represents $\mathbb P^1\times \mathbb P^1$. 
However, when using such figures to describe Calabi--Yau threefolds, the same 
polygons represent something else. Indeed, a given planar vertex of coordinates drawn as $(a,b)\in \mathbb R^2$
is to be interpreted as the head of a vector in $\mathbb R^3$ starting at the origin and ending in $(a,b,1)$. 
This works out correctly because the toric threefold is Calabi--Yau precisely when all the vectors defining its fan
end on the same hyperplane. The  most essential building blocks occurring in this work are:

\begin{center}
\begin{tikzpicture}[scale=.75]
\node at (0,3.5) {$ \{xy-zw=0\}=$ };
\draw (2,3) -- (3,3) -- (3,4) -- (2,4) -- (2,3);

\node at (3.6,3.5) { and };

\node at (6.3,3.5) {$ \{xy-z^2=0\}= $ };
\draw (8,3) -- (10,3) -- (9,4) -- (8,3);
\draw (8,3) -- (9,4);

\end{tikzpicture}
\end{center}
\vspace{-2cm}

With sample resolutions (= subdivisions into triangles of area $1/2$):

\begin{center}
\begin{tikzpicture}[scale=.75]
\node at (-1.2,3.5) {$ \Tot(\mathcal O_{\mathbb P^1}(-1)\oplus\mathcal O_{\mathbb P^1}(-1)) =$ };
\draw (2,3) -- (3,3) -- (3,4) -- (2,4) -- (2,3);
\draw[green](2,3)--(3,4);
\node at (7.3,3.5) {$ \Tot(\mathcal O_{\mathbb P^1}(-2)\oplus\mathcal O_{\mathbb P^1}) = $ };
\draw (10,3) -- (12,3) -- (11,4) -- (10,3);
\draw (10,3) -- (11,4);
\draw[red] (11,3) -- (11,4);
\end{tikzpicture}
\end{center}
\vspace{-2cm}

$\mathbb C^3/\mathbb Z_3$ singularity  \(\frac{1}{3}(1,1,1)\) \quad  and \quad  $\mathbb C^3/\mathbb Z_4$ singularity \(\frac{1}{4}(1, 1, 2)\)\\
  
With sample resolutions:

$ \Tot(\mathcal O_{\mathbb P^1}(-3)\oplus\mathcal O_{\mathbb P^1}(1)) = $
\begin{tikzpicture}[scale=0.75]
\draw		(0,0) -- (1,0);
\draw		(0,0) -- (0,-1);
\draw		(0,-1) -- (-1,1);
\draw		(1,0) -- (-1,1);
\draw[cyan]	(0,0) -- (-1,1);
\draw[dashed]	(0,-1) -- (1,0);
\end{tikzpicture}
$ \Tot(\mathcal O_{\mathbb P^1}(-4)\oplus\mathcal O_{\mathbb P^1}(2)) = $
\begin{tikzpicture}[scale=0.75]
\draw (0,0) -- (0,-1);
\draw (0,0) -- (2,-1);
\draw (0,-1) -- (-1,1);
\draw (2,-1) -- (-1,1);
\draw[cyan]	(0,0) -- (-1,1);
\draw[dashed]	(0,-1) -- (2,-1);
\end{tikzpicture}

There is some freedom in drawing the polygon, but one ought to first give 
the underlying lattice. In this work we will always use integer lattices. Thus, planar
polygons differing by an $GL(2,\mathbb Z)$ transformation determine the same complex toric variety, 
and 3-dimensional fans differing by an $GL(3,\mathbb Z)$ transformation determine the same toric variety.

\begin{example} Observe that the Calabi--Yau threefolds given by the following 2 polytopes are isomorphic.

\begin{center}
\begin{tikzpicture}
\draw (-1,1) -- (-1,0) -- (1,0) -- cycle;
\draw[red] (0,0) -- (-1,1);
\end{tikzpicture}
\quad and \quad
\begin{tikzpicture}
\draw (-1,0) -- (1,0) -- (0,1) -- cycle;
\draw[red] (0,0) -- (0,1);
\end{tikzpicture}
\end{center}

The rays of their fans are given by vectors in $N = \mathbb{Z}^3$ as:
\begin{itemize}
    \item ($\Sigma_1$): $v_1=(0,0,1)$, $v_2=(1,0,1)$, $v_3=(2,0,1)$, $v_4=(0,1,1)$.\\
    $$v_1-2 v_2+v_3 =0$$    
    \item ($\Sigma_2$): $w_1=(-1,0,1)$, $w_2=(0,0,1)$, $w_3=(1,0,1)$, $w_4=(0,1,1)$.
    $$w_1-2w_2+w_3=0$$
\end{itemize}
The matrix must do $A(v_i) =w_i$

Define a linear transformation $\phi: \mathbb{Z}^3 \to \mathbb{Z}^3$ by the matrix $A$:
\[
A = \begin{pmatrix} 1 & 1 & -1 \\ 0 & 1 & 0 \\ 0 & 0 & 1 \end{pmatrix} \text{.}
\]

Then $A\in  \text{SL}(3, \mathbb{Z})$, and  transforms the rays of $\Sigma_1$ into the rays of $\Sigma_2$. Specifically, 
$A(\Sigma_1) = \Sigma_2$ as sets of rays.

Since there exists a lattice automorphism $\phi$ that maps the set of rays of $\Sigma_1$ 
to the set of rays of $\Sigma_2$, and this transformation preserves the simplicial cone structure, the given toric varieties  are isomorphic.

\end{example}

We now provide an example of how to obtain the Calabi--Yau threefold from the polygon.

\begin{example}\label{strip24}
Consider the polygon 
\begin{center}\begin{tikzpicture}[scale=.75]
\foreach \x in {2,3,4,5,6} {
  \foreach \y in {3,4} {
    \node at (\x,\y) [circle,fill=black,scale=.1] {};
  }
}
\draw (6,3) -- (4,4) -- (2,4) -- (2,3) -- (6,3);
\end{tikzpicture}\end{center}
\vspace{-2cm}

The given polygon is determined by the vertices $(0,0),  (4,0), (0,1), (2,1). $
By definition, the vectors defining the corresponding threefold  $C_{2,4}$ are:
$$v_1= (0,0,1), v_2= (4,0,1), v_3= (0,1,1), v_4=(2,1,1). $$
Hence, these 4 vectors define a fan representing a Calabi--Yau threefold whose equation we wish 
to obtain. 

By the basic rules of toric geometry, one can obtain such an equation by 
finding the dual polytope. Such dual polytope is given  by the inward normals to the faces of the fan. 
So, we proceed to calculate inward normals. 
For the face $F_{1,2}$ , we have the normal vector  
$$n_{1,2}= \det\begin{pmatrix} 
i& j& k\\
0& 0& 1\\
4& 0 & 1\end{pmatrix}=4j,\quad \textnormal{but we take the shorter vector}\quad   A=j,$$
observing that we need the inward normals, taken with the smallest length that makes such 
a vector end at a lattice point inside the fan.

Next, for the face $F_{1,3}$ , we have the normal vector given by 
$$n_{1,3}= \det\begin{pmatrix} 
i& j& k\\
0& 0& 1\\
0& 1 & 1\end{pmatrix}=-i, \text{the inward pointing vector is} +i, \quad B= i$$
 for the face $F_{2,4}$ , we have the normal vector given by 
$$n_{2,4}= \det\begin{pmatrix} 
i& j& k\\
4& 0& 1\\
2& 1 & 1\end{pmatrix}=-i-2j+4k, \quad \text{thus} \quad C= -i-2j+4k,$$
and
 for the face $F_{3,4}$ , we have the normal vector given by 
$$n_{3,4}= \det\begin{pmatrix} 
i& j& k\\
0& 1& 1\\
2& 1 & 1\end{pmatrix}=2j-2k, \text{inner is} -2j+2k, \quad D= -2j+2k.$$

We now explore the relations among these vectors, expecting to  obtain 
a single equation for the fan (so as to determine a 3-dimensional variety 
out of the 4 variables). We have:
$B+C=2 A+4D,$ and taking exponentials named $x= e^B, y=e^C, z=e^A, w=e^D$ we get

$$xy= e^{B+C} = e^{2A+4D} = z^2w^4.$$
We conclude that the polygon with vertices 
$(0,0),  (4,0), (0,1), (2,1) $ represents the singular Calabi--Yau threefold
$$\boxed{xy-z^2w^4=0}$$

\end{example}

\section{Toy Example: a Triangle with One Internal Vertex}\label{toy}

We will call the following figure $T_{2,4}$. It appears as F13 in the listed favourite polygons of Section \ref{favpolys}.

\begin{center}
\begin{tikzpicture}[scale=0.75]
\foreach \x in {0,1,2,3,4} {
  \foreach \y in {0,1,2} {
    \node at (\x,\y) [circle,fill=black,scale=.1] {};
  }
}
\draw (0,0) -- (4,0) -- (0,2) -- cycle;
\end{tikzpicture}
\end{center}

This figure represents a Calabi--Yau threefold that has a quotient singularity 
obtained by the action of a group of order 8, in fact, by $\mathbb Z_2\times \mathbb Z_4$.
We consider first the auxiliary figure of the strip $C_{2,4}$ discussed in Example \ref{strip24}

\begin{center}
\begin{tikzpicture}[scale=0.75]
\foreach \x in {0,1,2,3,4} {
  \foreach \y in {0,1,2} {
    \node at (\x,\y) [circle,fill=black,scale=.1] {};
  }
}
\draw (0,0) -- (4,0) -- (0,2) -- cycle;
\draw (0,1) -- (2,1);
\draw[red] (1,1) -- (0,2);
\end{tikzpicture}
\end{center}

The triangulations of $C_{2,4}$ were given in \cite{GKMR}; they are:\\

\begin{tikzpicture}[scale=.75]
\draw[red] (1,0) -- (0,1); 
\draw[red] (2,0) -- (0,1);
\draw[red] (3,0) -- (0,1);
\draw[green] (4,0) -- (0,1);
\draw[red] (4,0) -- (1,1);
\draw[black] (0,0) -- (4,0) -- (2,1) -- (0,1) --cycle; 
\end{tikzpicture}
\begin{tikzpicture}[scale=.75]
\draw[red] (1,0) -- (0,1);
\draw[red] (2,0) -- (0,1);
\draw[green] (3,0) -- (0,1);
\draw[green] (3,0) -- (1,1);
\draw[green] (4,0) -- (1,1);
\draw[black] (0,0) -- (4,0) -- (2,1) -- (0,1) --cycle; 
\end{tikzpicture}
\begin{tikzpicture}[scale=.75]
\draw[red] (1,0) -- (0,1);
\draw[red] (2,0) -- (0,1);
\draw[green] (3,0) -- (0,1);
\draw[red] (3,0) -- (1,1);
\draw[green] (3,0) -- (2,1);
\draw[black] (0,0) -- (4,0) -- (2,1) -- (0,1) --cycle; 
\end{tikzpicture} 
\linebreak

\begin{tikzpicture}[scale=.75]
\draw[red] (1,0) -- (0,1);
\draw[green] (2,0) -- (0,1);
\draw[green] (2,0) -- (1,1);
\draw[red] (3,0) -- (1,1);
\draw[green] (4,0) -- (1,1);
\draw[black] (0,0) -- (4,0) -- (2,1) -- (0,1) --cycle; 
\end{tikzpicture}
\begin{tikzpicture}[scale=.75]
\draw[red] (1,0) -- (0,1);
\draw[green] (2,0) -- (0,1);
\draw[green] (2,0) -- (1,1);
\draw[green] (3,0) -- (1,1);
\draw[green] (3,0) -- (2,1);
\draw[black] (0,0) -- (4,0) -- (2,1) -- (0,1) --cycle; 
\end{tikzpicture}
\begin{tikzpicture}[scale=.75]
\draw[red] (1,0) -- (0,1);
\draw[green] (2,0) -- (0,1);
\draw[red] (2,0) -- (1,1);
\draw[green] (2,0) -- (2,1);
\draw[red] (3,0) -- (2,1);
\draw[black] (0,0) -- (4,0) -- (2,1) -- (0,1) --cycle; 
\end{tikzpicture} 
\linebreak

\begin{tikzpicture}[scale=.75]
\draw[green] (1,0) -- (0,1);
\draw[green] (1,0) -- (1,1);
\draw[red] (2,0) -- (1,1);
\draw[red] (3,0) -- (1,1);
\draw[green] (4,0) -- (1,1);
\draw[black] (0,0) -- (4,0) -- (2,1) -- (0,1) --cycle; 
\end{tikzpicture}
\begin{tikzpicture}[scale=.75]
\draw[green] (1,0) -- (0,1);
\draw[green] (1,0) -- (1,1);
\draw[red] (2,0) -- (1,1);
\draw[green] (3,0) -- (1,1);
\draw[green] (3,0) -- (2,1);
\draw[black] (0,0) -- (4,0) -- (2,1) -- (0,1) --cycle; 
\end{tikzpicture}
\begin{tikzpicture}[scale=.75]
\draw[green] (1,0) -- (0,1);
\draw[green] (1,0) -- (1,1);
\draw[green] (2,0) -- (1,1);
\draw[green] (2,0) -- (2,1);
\draw[red] (3,0) -- (2,1);
\draw[black] (0,0) -- (4,0) -- (2,1) -- (0,1) --cycle; 
\end{tikzpicture} 
\linebreak

\begin{tikzpicture}[scale=.75]
\draw[green] (1,0) -- (0,1);
\draw[red] (1,0) -- (1,1);
\draw[green] (1,0) -- (2,1);
\draw[red] (2,0) -- (2,1);
\draw[red] (3,0) -- (2,1);
\draw[black] (0,0) -- (4,0) -- (2,1) -- (0,1) --cycle; 
\end{tikzpicture}
\begin{tikzpicture}[scale=.75]
\draw[green] (0,0) -- (1,1);
\draw[red] (1,0) -- (1,1);
\draw[red] (2,0) -- (1,1);
\draw[red] (3,0) -- (1,1);
\draw[green] (4,0) -- (1,1);
\draw[black] (0,0) -- (4,0) -- (2,1) -- (0,1) --cycle; 
\end{tikzpicture}
\begin{tikzpicture}[scale=.75]
\draw[green] (0,0) -- (1,1);
\draw[red] (1,0) -- (1,1);
\draw[red] (2,0) -- (1,1);
\draw[green] (3,0) -- (1,1);
\draw[green] (3,0) -- (2,1);
\draw[black] (0,0) -- (4,0) -- (2,1) -- (0,1) --cycle; 
\end{tikzpicture} 
\linebreak

\begin{tikzpicture}[scale=.75]
\draw[black] (0,0) -- (4,0) -- (2,1) -- (0,1) --cycle; 
\draw[green] (0,0) -- (1,1);
\draw[red] (1,0) -- (1,1);
\draw[green] (2,0) -- (1,1);
\draw[green] (2,0) -- (2,1);
\draw[red] (3,0) -- (2,1);
\end{tikzpicture}
\begin{tikzpicture}[scale=.75]
\draw[black] (0,0) -- (4,0) -- (2,1) -- (0,1) --cycle; 
\draw[green] (0,0) -- (1,1);
\draw[green] (1,0) -- (1,1);
\draw[green] (1,0) -- (2,1);
\draw[red] (2,0) -- (2,1);
\draw[red] (3,0) -- (2,1);
\end{tikzpicture}
\begin{tikzpicture}[scale=.75]
\draw[black] (0,0) -- (4,0) -- (2,1) -- (0,1) --cycle; 
\draw[red] (0,0) -- (1,1);
\draw[green] (0,0) -- (2,1);
\draw[red] (1,0) -- (2,1);
\draw[red] (2,0) -- (2,1);
\draw[red] (3,0) -- (2,1);
\end{tikzpicture} 
\linebreak

The following  triangulations of $T_{2,4}$  do not come from those of $C_{2,4}$.

\vspace{12pt}

\begin{tikzpicture}[scale=0.7]
\draw[red] (0,1) -- (1,0);
\draw[red] (1,0) -- (1,1);
\draw[red] (1,1) -- (2,0);
\draw[red] (1,1) -- (3,0);
\draw[red] (1,1) -- (2,1);
\draw[green] (1,0) -- (0,2);
\draw[green] (1,1) -- (4,0);
\draw[cyan] (1,1) -- (0,2);
\draw (0,0) -- (4,0) -- (0,2) -- cycle; 
\node[xshift=-7pt, yshift=-6pt] at (4,2) {{\bf T1}};
\end{tikzpicture}
\hfill
\begin{tikzpicture}[scale=0.7]
\draw[red] (0,1) -- (1,0);
\draw[red] (1,0) -- (1,1);
\draw[red] (1,1) -- (2,1);
\draw[red] (2,1) -- (3,0);
\draw[green] (1,0) -- (0,2);
\draw[green] (1,1) -- (2,0);
\draw[green] (2,0) -- (2,1);
\draw[cyan] (1,1) -- (0,2);

\draw (0,0) -- (4,0) -- (0,2) -- cycle; 
\node[xshift=-7pt, yshift=-6pt] at (4,2) {{\bf T2}};
\end{tikzpicture}
\hfill
\begin{tikzpicture}[scale=0.7]
\draw[red] (0,1) -- (1,0);
\draw[red] (2,0) -- (2,1);
\draw[red] (2,1) -- (3,0);
\draw[green] (1,0) -- (0,2);
\draw[green] (1,0) -- (2,1);
\draw[cyan] (1,1) -- (0,2);
\draw[cyan] (1,1) -- (1,0);
\draw[cyan] (1,1) -- (2,1);
\draw (0,0) -- (4,0) -- (0,2) -- cycle; 
\node[xshift=-7pt, yshift=-6pt] at (4,2) {{\bf T3}};
\end{tikzpicture}

\vspace{12pt}

\begin{tikzpicture}[scale=0.7]
\draw[red] (1,1) -- (2,0);
\draw[red] (1,1) -- (3,0);
\draw[red] (3,0) -- (2,1);
\draw[green] (1,0) -- (1,1);
\draw[green] (1,1) -- (0,1);
\draw[green] (1,0) -- (0,1);
\draw[green] (0,2) -- (3,0);
\draw[cyan] (1,1) -- (0,2);
\draw (0,0) -- (4,0) -- (0,2) -- cycle; 
\node[xshift=-7pt, yshift=-6pt] at (4,2) {{\bf T4}};
\end{tikzpicture}
\hfill
\begin{tikzpicture}[scale=0.7]
\draw[red] (0,1) -- (1,1);
\draw[red] (1,0) -- (1,1);
\draw[red] (1,1) -- (2,0);
\draw[red] (1,1) -- (3,0);
\draw[red] (3,0) -- (2,1);
\draw[green] (0,0) -- (1,1);
\draw[green] (3,0) -- (0,2);
\draw[cyan] (1,1) -- (0,2);
\draw (0,0) -- (4,0) -- (0,2) -- cycle; 
\node[xshift=-7pt, yshift=-6pt] at (4,2) {{\bf T5}};
\end{tikzpicture}
\hfill
\begin{tikzpicture}[scale=0.7]
\draw[red] (0,1) -- (1,0);
\draw[red] (1,0) -- (1,1);
\draw[red] (1,1) -- (2,0);
\draw[red] (1,1) -- (3,0);
\draw[red] (3,0) --(2,1);
\draw[green] (1,0) -- (0,2);
\draw[green] (3,0) -- (0,2);
\draw[cyan] (1,1) -- (0,2);
\draw (0,0) -- (4,0) -- (0,2) -- cycle; 
\node[xshift=-7pt, yshift=-6pt] at (4,2) {{\bf T6}};
\end{tikzpicture}

\vspace{12pt}

\begin{tikzpicture}[scale=0.7]
\draw[red] (1,0) -- (0,1);
\draw[red] (0,1) -- (1,1);
\draw[red] (1,1) -- (3,0);
\draw[red] (3,0) -- (2,1);
\draw[green] (0,1) -- (2,0);
\draw[green] (2,0) -- (1,1);
\draw[green] (3,0) -- (0,2);
\draw[cyan] (1,1) -- (0,2);
\draw (0,0) -- (4,0) -- (0,2) -- cycle; 
\node[xshift=-7pt, yshift=-6pt] at (4,2) {{\bf T7}};
\end{tikzpicture}
\hfill
\begin{tikzpicture}[scale=0.7]
\draw[red] (0,1) -- (1,0);
\draw[red] (1,0) -- (1,1);
\draw[red] (1,1) -- (2,0);
\draw[green] (0,2) -- (1,0);
\draw[green] (1,1) -- (3,0);
\draw[green] (1,1) -- (2,1);
\draw[green] (3,0) -- (2,1);
\draw[cyan] (1,1) -- (0,2);
\draw (0,0) -- (4,0) -- (0,2) -- cycle; 
\node[xshift=-7pt, yshift=-6pt] at (4,2) {{\bf T8}};
\end{tikzpicture}
\hfill
\begin{tikzpicture}[scale=0.7]
\draw[red] (0,1) -- (1,0);
\draw[red] (0,1) -- (2,0);
\draw[red] (3,0) -- (2,1);
\draw[green] (0,1) -- (3,0);
\draw[green] (0,2) -- (3,0);
\draw[cyan] (1,1) -- (0,2);
\draw[cyan] (1,1) -- (0,1);
\draw[cyan] (1,1) -- (3,0);
\draw (0,0) -- (4,0) -- (0,2) -- cycle; 
\node[xshift=-7pt, yshift=-6pt] at (4,2) {{\bf T9}};
\end{tikzpicture}

Observe that $T1, \cdots,  T9$ all contain  edges coloured blue, which 
are neither inside triangles nor inside parallelograms.

\subsection{Resolutions and Blowups }

Our goal now is to use blowups as a tool to help us decide which 
resolutions  of $T_{2,4}$  are most likely to be the initial configuration 
from which the singularity itself  originated. Using physics terminology we may
informally think of such resolution  as the most probable past of the singularity. 

Let us first recall some basic properties of Chern classes. 
It is enough for us to discuss the case of surfaces, and think of 
the Chern classes of the compact surfaces represented by the planar  polytopes. 

The total Chern class (\(c = 1 + c_1 + c_2\)) of the blown-up surface \(\tilde {Y}\)
is calculated using the formulas below, where \(p\colon \tilde{Y} \to Y\) is the projection map, 
\(E\) is the exceptional divisor, and \(E^2 = -1\).

\[c(\tilde{Y})=(1+p^{*}c_{1}(Y)-E)(1+p^{*}c_{2}(Y)+E^{2}).\]

Therefore, for the first Chern Class:
$c_{1}(\tilde{Y})=p^{*}c_{1}(Y)-E,$
and for the second Chern Class:
\[c_{2}(\tilde{Y})=p^{*}c_{2}(Y)+E^{2}=c_{2}(Y)+1.\]
Since \(c_2(Y)\) is the Euler characteristic \(\chi(Y)\) for surfaces, blowing up at one point increases the second Chern class by 1.
Using physics jargon, we argue that this increase in second Chern class requires { energy} to 
be accomplished. Therefore, we will look at those resolutions appearing after more blowups as 
being less likely to occur. 

Note that we also have that the Chern classes of $X = \Tot(\omega_Y)$ are 
$c_1(X) = c_3(X) =0$ and $c_2(X)  =  -c_1^2(Y)+c_2(Y)$ (see \cite[Lem.\thinspace 3.4]{BGGH}).
Hence, blowing up the surface $Y$ also increases the second Chern class of the threefold $X$,
so that from the threefold viewpoint, the resolutions appearing after more blowups are 
less likely as well. 

The blowups of  $T_{2,4}$  triangulations  are described on  the table below.
\newcommand{\tablescale}{0.7}
\begin{table}
\begin{tblr}{Q[c,b]|Q[c,m]|Q[c,m]}
\hline
Triangulation & Blowups & Polytope \\
\hline 
$T_0$ & 0 & 
\begin{tikzpicture}[scale=\tablescale]
\foreach \x in {-4,-3,-2,-1,0} {
  \foreach \y in {0,1,2} {
  }
}
\draw (0,0) -- (-2,1) -- (-4,2) -- (-4,0) -- (0,0);
\draw[red]  (-4,2) -- (-3,1);
\draw[green]  (-4,0) -- (-3,1);
\draw[red]  (-3,1) -- (-3,0);
\draw[red]  (-4,1) -- (-3,1);
\draw[red]  (-3,1) -- (-2,1);
\draw[red]  (-3,1) -- (-2,0);
\draw[red]  (-3,1) -- (-1,0);
\draw[green]  (-3,1) -- (0,0);
\end{tikzpicture}\\
\hline
$\multitilde{2}{T_1}$ & 2 & 
\begin{tikzpicture}[scale=\tablescale]
\foreach \x in {-4,-3,-2,-1,0} {
  \foreach \y in {0, 1, 2} {
  }
}
\draw (-4,2) --(0,0) -- (-2,1)  ;
\draw (-3,0) -- (0,0);
\draw[green] (-4,2)--(-3,0);
\draw[red] (-3,1)--(-3,0);
\draw[cyan] (-4,2)--(-3,1);
\draw[red](-3,1)--(-2,1);
\draw[red](-3,1)--(-2,0);
\draw[red](-3,1)--(-1,0);
\draw[green](-3,1)--(0,0);
\end{tikzpicture} \\
\hline
$\multitilde{2}{T_5}$ & 2 & 
\begin{tikzpicture}[scale=\tablescale]
\foreach \x in {2,3,4,5,6} {
  \foreach \y in {0,1,2} {
  }  
}
\draw   (2,2) -- (2,0) -- (5,0);
\draw[green]  (2,2) -- (5,0);
\draw[green]  (2,0) -- (3,1);
\draw[red]  (3,1) -- (3,0);
\draw[cyan] (2,2) -- (3,1);
\draw[red]  (2,1) -- (3,1);
\draw[red]  (3,1) -- (4,0);
\draw[red]  (3,1) -- (5,0);
\end{tikzpicture}\\
\hline
$\multitilde{3}{T_4}$ & 3 & 
\begin{tikzpicture}[scale=\tablescale]
\foreach \x in {-4,-3,-2,-1,0} {
  \foreach \y in {0,1,2} {
  } 
}
\draw   (-4,2)--(-4,1);
\draw  (-3,0) -- (-1,0);
\draw[green]  (-4,2) -- (-1,0);
\draw[green]  (-4,1) -- (-3,0);
\draw[green]  (-3,1) -- (-3,0);
\draw[cyan]   (-4,2) -- (-3,1);
\draw[green]  (-4,1) -- (-3,1);
\draw[red]  (-3,1) -- (-2,0);
\draw[red]  (-3,1) -- (-1,0);
\end{tikzpicture}\\
\hline
$\multitilde{3}{T_8}$ & 3 & 
\begin{tikzpicture}[scale=\tablescale]
\foreach \x in {2,3,4,5,6} {
  \foreach \y in {0,1, 2} {
  } 
}
\draw  (4,1) -- (2,2);
\draw  (3,0) -- (5,0);
\draw[green]  (2,2) -- (3,0);
\draw[red]  (3,1) -- (3,0);
\draw[cyan] (2,2) -- (3,1);
\draw[green]  (3,1) -- (4,1);
\draw[red]  (3,1) -- (4,0);
\draw[green]  (5,0) -- (3,1);
\draw[green]  (4,1) -- (5,0);
\end{tikzpicture}\\
\hline
$\multitilde{4}{T_2}$ & 4 &
\begin{tikzpicture}[scale=\tablescale]
\foreach \x in {2,3,4,5,6} {
  \foreach \y in {0,1,2} {
  }
}
\draw  (4,1) -- (2,2);
\draw  (3,0) -- (4,0);
\draw[green]  (2,2) -- (3,0);
\draw[red]  (3,1) -- (3,0);
\draw[cyan] (2,2) -- (3,1);
\draw[red]  (3,1) -- (4,1);
\draw[green]  (3,1) -- (4,0);
\draw[green]  (4,1) -- (4,0);
\end{tikzpicture} \\
\hline
$\multitilde{4}{T_6}$ & 4 &
\begin{tikzpicture}[scale=\tablescale]
\foreach \x in {8,9,10,11,12} {
  \foreach \y in {0,1,2} {
  }
}
\draw (9,0) -- (11,0);
\draw[green]  (8,2) -- (11,0);
\draw[red]  (9,1) -- (9,0);
\draw[cyan] (8,2) -- (9,1);
\draw[green]  (8,2) -- (9,0);
\draw[red]  (9,1) -- (10,0);
\draw[red]  (9,1) -- (11,0);
\end{tikzpicture} \\
\hline
$\multitilde{4}{T_7}$ & 4 &
\begin{tikzpicture}[scale=\tablescale]
\foreach \x in {-4,-3,-2,-1,0} {
  \foreach \y in {0,1,2} {
  }
}
\draw (-4,2) -- (-4,1);
\draw (-2,0) -- (-1,0);
\draw[green]  (-4,2) -- (-1,0);
\draw[green]  (-4,1) -- (-2,0);
\draw[cyan] (-4,2) -- (-3,1);
\draw[red]  (-4,1) -- (-3,1);
\draw[green]  (-3,1) -- (-2,0);
\draw[red]  (-3,1) -- (-1,0);
\end{tikzpicture}\\
\hline
$\multitilde{5}{T_3}$ & 5 & 
\begin{tikzpicture}[scale=\tablescale]
\foreach \x in {8,9,10,11,12} {
  \foreach \y in {0,1,2} {
  }
}
\draw (10,1) -- (8,2);
\draw[green]  (8,2) -- (9,0);
\draw[cyan] (9,1) -- (9,0);
\draw[cyan] (8,2) -- (9,1);
\draw[cyan] (9,1) -- (10,1);
\draw[green]  (9,0) -- (10,1);
\end{tikzpicture}\\
\hline
 $\multitilde{5}{T_9}$ & 5 & 
\begin{tikzpicture}[scale=\tablescale]
\foreach \x in {8,9,10,11,12} {
  \foreach \y in {0,1,2} {
  }
}
\draw (8,2) -- (8,1);
\draw[green]  (8,2) -- (11,0);
\draw[cyan] (8,1) -- (9,1);
\draw[cyan] (8,2) -- (9,1);
\draw[green]  (8,1) -- (11,0);
\draw[cyan] (9,1) -- (11,0);
\end{tikzpicture}\\
\hline
\end{tblr}
\caption*{Triangulations of $T_{2,4}$ with  corners chopped off}
\end{table}

Thus,  we see here that  the probabilities are higher according to 
  the larger numbers of red lines. We conclude that, the  most likely configuration to occur is
  actually  $T_0$ which has 6 red lines. We observe that this conclusion is in agreement with 
the one obtained in \cite{GSTV}. While  different colour choices were used in that reference, the 
geometric conclusion coincides, reinforcing the observation  that $(-2,0)$ are the most likely lines to occur.

\clearpage
\noindent We finish this section with  observations about charts containing blue lines.\\

\noindent{\bf The case of  $(-3,1)$ lines}\\

\begin{tikzpicture}[scale=1][h]
\foreach \x in {-1,0,1} {
  \foreach \y in {-1,0,1} {
    \node at (\x,\y) [circle,fill=black,scale=.1] {};
  }
}
\draw  (0,-1) -- (1,0) -- (-1,1) -- cycle;
\draw[cyan] (0,0) -- (-1,1);
\draw[cyan] (0,0) -- (1,0);
\draw[cyan] (0,0) -- (0,-1);
\node[xshift=18pt, yshift=-6pt] at (2,1) 
{${\bf \multitilde{5}{T_3} \simeq F1 \simeq \multitilde{5}{T_9}}$} ;
\end{tikzpicture}

\vspace{12pt}

Without the internal point, this triangle gives just the 
quotient singularity $\mathbb C^3/\mathbb Z_3$ 
with the diagonal action, that is, with weights $(1,1,1)$. 
The resolution given by the radial triangulation is described 
by the fan corresponding to this picture. It has rays: 
$$v_1= (0,-1,1), v_2=  (1,0,1), v_3= (-1,1,1), v_4=(0,0,1)$$
these satisfy
$$v_1+ v_2+v_3 - 3v_4= 0$$
which means that the variety can be written as a quotient of $\mathbb C^4$ 
(minus the subvariety $x=y=z=0$ which does not form a subcone)
by 
the action 

$$t \cdot (v_1,v_2,v_3,v_4) = (tv_1,tv_2,tv_3,t^{-3}v_4)$$
without the $v_4$ vector coming from the internal point, the 
projectivisation of this action gives $\mathbb P^2$
and with the $v_4$ included it gives $\mathcal O_{\mathbb P^2}(-3)$ 
because of the weight $-3$ of the action on $v_4$. 
Note that  $\mathcal O_{\mathbb P^2}(-3) = K_{\mathbb P^2}$ 
the canonical bundle of  $\mathbb P^2$ (the only possibility for 
embedding $\mathbb P^2$ in a CY3). \\

\begin{tikzpicture}[scale=1]
\foreach \x in {-1,0,1,2} {
  \foreach \y in {-1,0,1} {
  \node at (\x,\y) [circle,fill=black,scale=.1] {};
  }
}

\draw[cyan] (0,0) -- (-1,0);
\draw[cyan] (0,0) -- (-1,1);
\draw[cyan] (0,0) -- (2,-1);
\draw (-1,1) -- (-1,0) -- (2,-1) -- cycle;
\node[xshift=0, yshift=-6pt] at (3,1) {${\bf \multitilde{5}{ T_9}}$} ; 
\end{tikzpicture}

In this case the vectors are
$$v_1= (-1,0,1), v_2=  (2,-1,1), v_3= (-1,1,1), v_4=(0,0,1),$$
and they satisfy 
$$v_1+v_2+v_3-3v_4=0 .$$
Therefore, 
the variety can be written as a quotient of $\mathbb C^4$ 
 (minus the subvariety which does not form a subcone, there are 3 subcones)
 by 
 the action 
  $$t \cdot (v_1,v_2,v_3,v_4) = (tv_1,tv_2,tv_3,t^{-3}v_4).$$
  The same as the one obtained in $\multitilde{5}{T_3}$.
  Indeed, let multiplication by
\(\tiny{\left[\begin{matrix}2&1&0\\ -1&0&0\\ 0&0&1\end{matrix}\right]}\)
   takes \{(0,-1,1), (1,0,1), (-1,1,1), (0,0,1)\}
  to \{(-1,0,1),  (2,-1,1),  (-1,1,1), (0,0,1)\}.\\

\noindent{\bf The case of $(-4,2)$ lines}:\label{sing42}\\

\begin{tikzpicture}[scale=0.7]
\foreach \x in {-1,0,1,2} {
  \foreach \y in {-1,0,1} {
    \node at (\x,\y) [circle,fill=black,scale=.1] {};
  }
}
\draw[cyan] (0,0) -- (-1,1);
\draw (-1,1) -- (0,-1) -- (0,0) -- (2,-1) -- cycle;
\node[xshift=0, yshift=-6pt] at (3,1) {${\bf \multitilde{4}{T_6}}$};
\end{tikzpicture}

Here the vectors of the fan are
$$v_1= (0,-1,1), v_2= (2,-1,1), v_3= (-1,1,1), v_4= (0,0,1),$$
and they satisfy the equation is 
$$v_1+v_2+2v_3-4 v_4=0.$$

Therefore, we get the action   $$t \cdot (v_1,v_2,v_3,v_4) = (tv_1,tv_2,t^2v_3,t^{-4}v_4),$$
which describes the singularity
\(\frac{1}{4}(1, 1, 2)\), obtained as the quotient space \(\mathbb{C}^3 / \mathbb{Z}_4\) under the group action:\((e^{2\pi i/4}\cdot x,\;e^{2\pi i/4}\cdot y,\;e^{2\pi i\cdot 2/4}\cdot z)=(ix,\;iy,\;-z)\). \\

Similar calculations show that, with the exception of $\multitilde{4}{T_6}$ which has a $(-4,2)$ line, 
all other  blue lines appearing in this section are of type $(-3,1)$

\section{16 Favourite Polygons}\label{favpolys}

Let us first define a couple of concepts that will be used in this section. \\

For a convex polytope \(P\) containing the origin in its interior, the {\it polar dual} \(P^{\circ }\) is defined as the set of all points \(y\) such that \(x \cdot y \leq 1\) for every \(x \in P\). \\

 We will also use the concept of number of moduli, but which we mean the number 
 of triangulations moduli isomorphism, or equivalently, the number of distinct CY3's depicted by such triangulations. 
  
  Two fans determine { isomorphic} toric varieties if and only if there is an isomorphism of the underlying lattices that maps one fan bijectively onto the other.
In terms of algebraic varieties, an isomorphism of the algebraic tori extends to an isomorphism of the corresponding toric varieties if and only if their defining fans are isomorphic under a linear automorphism of the lattice. 
More precisely:

\begin{definition}
Let $N_1$ and $N_2$ be two lattices of the same finite rank $n$. Let $\Delta_1 \subset N_{1} \otimes \mathbb{R}$ and $\Delta_2 \subset N_{2} \otimes \mathbb{R}$ be two fans. The induced toric varieties $X_{\Delta_1}$ and $X_{\Delta_2}$ are isomorphic if there exists a lattice isomorphism:
\[
\phi: N_1 \xrightarrow{\sim} N_2
\]

The linear extension of the lattice isomorphism $\phi_\mathbb{R}: N_1 \otimes \mathbb{R} \to N_2 \otimes \mathbb{R}$ must satisfy the following geometric conditions:
\begin{itemize}
    \item For every cone $\sigma_1 \in \Delta_1$, its image $\phi_\mathbb{R}(\sigma_1)$ must be a cone $\sigma_2 \in \Delta_2$.
    \item  The assignment $\sigma_1 \mapsto \phi_\mathbb{R}(\sigma_1)$ defines a bijection between the sets of cones $\Delta_1 \to \Delta_2$.
\end{itemize}
When such a $\phi$ exists, the fans $\Delta_1$ and $\Delta_2$ are said to be {\it isomorphic fans}.
\end{definition}

Thus, two fans determine isomorphic toric varieties if and only if they are isomorphic under a $\text{GL}(n, \mathbb{Z})$ transformation of the underlying lattices.\\

In the following list, all 2D reflexive polytopes containing one lattice point in its interior, they  are depicted up to $SL(2,\mathbb Z)$ transformations.  
Regarded together with the internal dashed edges, figures $F1, \dots, F16$  form toric diagrams of smooth complex surfaces, and they
represent Del Pezzo surfaces.  
\begin{notation}\label{polys} For $i=1, \dots, 16$, we denote by $Pi$ the polygon which forms the boundary of figure $Fi$ below (with 
its corresponding internal lattice point).
\end{notation}
The polygons $Pi$ are dual to $P(17-i)$  for $i=1,\dots,6$ 
 while $P7, P8, P9,$ and $ P10$ are self-dual, they appear in \cite{CKYZ, GP}.  Updownarrow means \lq\lq polar dual to\rq\rq.
\\

\noindent \begin{minipage}{.1\textwidth}
  \begin{tikzpicture}

  \draw \stB (0,0) -- (0,1);
  \draw \stB (0,0) -- (1,0);
  \draw \stB (0,0) -- (-1,-1);  

  \draw \stA (0,1) -- (1,0) -- (-1,-1) -- cycle;

  \end{tikzpicture}
\end{minipage}  \quad \quad\,
\begin{minipage}{.1\textwidth}
  \begin{tikzpicture}

  \draw \stB (0,0) -- (0,1);
  \draw \stB (0,0) -- (1,0);
  \draw \stB (0,0) -- (0,-1);
  \draw \stB (0,0) -- (-1,0);
  
  \draw \stA (0,1) -- (1,0) -- (0,-1) -- (-1,0) -- cycle;

  \end{tikzpicture}
\end{minipage}  \quad \quad \,
\begin{minipage}{.1\textwidth}
  \begin{tikzpicture}

  \draw \stB (0,0) -- (0,1);
  \draw \stB (0,0) -- (1,0);
  \draw \stB (0,0) -- (-1,0); 
  \draw \stB (0,0) -- (-1,-1);  

  \draw \stA (0,1) -- (1,0) -- (-1,-1) -- (-1,0) -- cycle;

  \end{tikzpicture}
\end{minipage}  \!\!\!
\begin{minipage}{.1\textwidth}
  \begin{tikzpicture}

  \draw \stB (0,0) -- (0,1);
  \draw \stB (0,0) -- (1,0);
  \draw \stB (0,0) -- (-1,0);
  \draw \stB (0,0) -- (-2,-1);  

  \draw \stA (0,1) -- (1,0) -- (-2,-1) -- (-1,0) -- cycle;

  \end{tikzpicture}
  
\end{minipage}  \quad\quad \quad \quad \quad
\begin{minipage}{.1\textwidth}
  \begin{tikzpicture}

  \draw \stB (0,0) -- (0,1);
  \draw \stB (0,0) -- (1,1);
  \draw \stB (0,0) -- (1,0);  
  \draw \stB (0,0) -- (0,-1); 
  \draw \stB (0,0) -- (-1,0); 
  \draw \stA (0,1) -- (1,1) -- (1,0) -- (0,-1) -- (-1,0) -- cycle;

  \end{tikzpicture}
  
\end{minipage}  \quad\quad \,\,
\begin{minipage}{.1\textwidth}
  \begin{tikzpicture}

  \draw \stB (0,0) -- (0,1);
  \draw \stB (0,0) -- (1,0);
  \draw \stB (0,0) -- (-1,-1);  
  \draw \stB (0,0) -- (-1,0); 
  \draw \stB (0,0) -- (-1,1); 
  \draw \stA (0,1) -- (1,0) -- (-1,-1) -- (-1,1) -- cycle;

  \end{tikzpicture}
\end{minipage}  
\hspace*{.4cm} F1\hspace{1.6cm} F2\hspace{1.6cm} F3\hspace{1.8cm} F4\hspace{1.6cm} F5\hspace{1.8cm} F6\\
\hspace*{.4cm} $\Updownarrow$ \hspace{1.8cm} $\Updownarrow$ \hspace{1.6cm} $\Updownarrow$ \hspace{1.9cm} 
$\Updownarrow$ \hspace{1.6cm} $\Updownarrow$ \hspace{1.8cm} $\Updownarrow$\\


\hspace*{-1cm}
\begin{minipage}{.1\textwidth}
  \begin{tikzpicture}

  \draw \stB (0,0) -- (0,1);
  \draw \stB (0,0) -- (1,0);
  \draw \stB (0,0) -- (2,-1);
  \draw \stB (0,0) -- (1,-1);
  \draw \stB (0,0) -- (0,-1);
  \draw \stB (0,0) -- (-1,-1);
  \draw \stB (0,0) -- (-1,0);
  \draw \stB (0,0) -- (-1,1);
  \draw \stB (0,0) -- (-1,2);
  \draw \stA (2,-1) -- (-1,-1) -- (-1,2) -- cycle;

  \end{tikzpicture}
  \end{minipage}\quad\quad\quad\quad\,\,
  \begin{minipage}{.1\textwidth}
    \begin{tikzpicture}

  \draw \stB (0,0) -- (0,1);
  \draw \stB (0,0) -- (1,1);
  \draw \stB (0,0) -- (1,0);
  \draw \stB (0,0) -- (1,-1);
  \draw \stB (0,0) -- (0,-1);
  \draw \stB (0,0) -- (-1,-1);
  \draw \stB (0,0) -- (-1,0);
  \draw \stB (0,0) -- (-1,1);

  \draw \stA (1,1) -- (1,-1) -- (-1,-1) -- (-1,1) -- cycle;
  \end{tikzpicture}
\end{minipage}  \quad \quad\,
\begin{minipage}{.1\textwidth}
  \begin{tikzpicture}
  \draw \stB (0,0) -- (0,1);
  \draw \stB (0,0) -- (1,0);
  \draw \stB (0,0) -- (1,-1);
  \draw \stB (0,0) -- (0,-1);
  \draw \stB (0,0) -- (-1,-1);
  \draw \stB (0,0) -- (-1,0);
  \draw \stB (0,0) -- (-1,1);
  \draw \stB (0,0) -- (-1,2);
  \draw \stA (1,0) -- (1,-1) -- (-1,-1) -- (-1,2) -- cycle;
  \end{tikzpicture}
\end{minipage}  \quad \quad\,
\begin{minipage}{.1\textwidth}
  \begin{tikzpicture}
  \draw \stB (0,0) -- (0,1);
  \draw \stB (0,0) -- (1,-1);
  \draw \stB (0,0) -- (0,-1);
  \draw \stB (0,0) -- (-1,-1);
  \draw \stB (0,0) -- (-1,0);
  \draw \stB (0,0) -- (-1,1);
  \draw \stB (0,0) -- (-1,2);
  \draw \stB (0,0) -- (-1,3);

  \draw \stA (1,-1) -- (-1,-1) -- (-1,3) -- cycle;
  \end{tikzpicture}
\end{minipage}  \quad \quad\,
\begin{minipage}{.1\textwidth}
  \begin{tikzpicture}
  \draw \stB (0,0) -- (0,1);
  \draw \stB (0,0) -- (1,0);
  \draw \stB (0,0) -- (1,-1); 
  \draw \stB (0,0) -- (0,-1);
  \draw \stB (0,0) -- (-1,-1);
  \draw \stB (0,0) -- (-1,0); 
  \draw \stB (0,0) -- (-1,1);
  \draw \stA (0,1) -- (1,0) -- (1,-1) -- (-1,-1) -- (-1,1) -- cycle;

  \end{tikzpicture}
\end{minipage}  \quad \quad\,
\begin{minipage}{.1\textwidth}
  \begin{tikzpicture}

  \draw \stB (0,0) -- (0,1);
  \draw \stB (0,0) -- (1,0);
  \draw \stB (0,0) -- (0,-1);
  \draw \stB (0,0) -- (-1,-1);  
  \draw \stB (0,0) -- (-1,0);
  \draw \stB (0,0) -- (-1,1);
  \draw \stB (0,0) -- (-1,2); 
  \draw \stA (1,0) -- (0,-1) -- (-1,-1) -- (-1,2) -- cycle;

  \end{tikzpicture}
\end{minipage}  \quad \quad\,

\hspace*{.3cm} F16\hspace{1.4cm} F15\hspace{1.4cm} F14\hspace{1.4cm} F13\hspace{1.6cm} F12\hspace{1.8cm} F11\\
  

\begin{minipage}{.1\textwidth}
  \begin{tikzpicture}
  \draw \stB (0,0) -- (0,1);
  \draw \stB (0,0) -- (1,1);
  \draw \stB (0,0) -- (1,0);
  \draw \stB (0,0) -- (0,-1);
  \draw \stB (0,0) -- (-1,-1);
  \draw \stB (0,0) -- (-1,0);

  \draw \stA (0,1) -- (1,1) -- (1,0) -- (0,-1) -- (-1,-1) -- (-1,0) -- cycle;

  \end{tikzpicture}
\end{minipage}  \quad \quad\,
\begin{minipage}{.1\textwidth}
  \begin{tikzpicture}

  \draw \stB (0,0) -- (-1,1);
  \draw \stB (0,0) -- (1,0);
  \draw \stB (0,0) -- (1,-1);
  \draw \stB (0,0) -- (0,-1);
  \draw \stB (0,0) -- (-1,-1);
  \draw \stB (0,0) -- (-1,0); 

  \draw \stA (-1,1) -- (1,0) -- (1,-1) --(-1,-1) -- cycle;

  \end{tikzpicture}
\end{minipage}  \quad \quad\,
\begin{minipage}{.1\textwidth}
  \begin{tikzpicture}

  \draw \stB (0,0) -- (0,1);
  \draw \stB (0,0) -- (1,1);
  \draw \stB (0,0) -- (1,0);  
  \draw \stB (0,0) -- (0,-1);
  \draw \stB (0,0) -- (-1,0);
  \draw \stB (0,0) -- (-1,1); 

  \draw \stA (1,1) -- (1,0) -- (0,-1) -- (-1,0) -- (-1,1) -- cycle;

  \end{tikzpicture}
\end{minipage}  \quad \quad\,
\begin{minipage}{.1\textwidth}
  \begin{tikzpicture}

  \draw \stB (0,0) -- (0,1);
  \draw \stB (0,0) -- (1,0);
  \draw \stB (0,0) -- (-1,-1);  
  \draw \stB (0,0) -- (-1,0);
  \draw \stB (0,0) -- (-1,1);
  \draw \stB (0,0) -- (-1,2); 

  \draw \stA (1,0) -- (-1,-1) -- (-1,2) -- cycle;

  \end{tikzpicture}
\end{minipage}  

\hspace{.1cm} F7\hspace{1.6cm} F8\hspace{1.7cm} F9\hspace{1.9cm} F10\\

\begin{example}\label{trigmod} For these 16 cases, we list the number of  
triangulations produced by our algorithm,  number of moduli and check for  convexity of a universal quilt: 

\begin{figure}[h]\label{trimod}
\begin{tabular} {c|c|c|c|c|c|c|c|c|c|c|c|c|c|c|c|c}\label{moduli}
Polygon & 1&2&3&4&5&6&7&8&9&10&11&12&13&14&15&16\\
\hline
 triangulations &1& 1& 2& 1& 5& 4& 18& 8& 13& 5& 18& 33& 24& 53& 64&79\\
 \hline
moduli  &1& 1& 2 & 1& 4& 4& 6 &8 &8 & 5&18 &19 &24 & 53 & 14&43\\
\hline
convex quilt &\checkmark &\checkmark&\checkmark&\checkmark  & & & &\checkmark&& \checkmark & & &\checkmark&\checkmark
&\checkmark&\checkmark
  \end{tabular}
\end{figure}
\end{example}  

\newpage

\subsection{Comparing Toric Threefolds Depicted by Triangles}

Geometric properties of a toric threefold $X_{\Sigma}$ are determined by the lattice vectors $v_i$ defining its fan $\Sigma$. 
For a 3-dimensional cone $\sigma = \text{cone}(v_1, v_2, v_3)$, the local patch is isomorphic to the quotient
 $\mathbb{C}^3 / G$, where the order of the group $G$ is given by the absolute determinant:
\[
|G| = | \det(v_1, v_2, v_3) |
\]

\begin{itemize}
    \item \textbf{Polygon P1:} 
     \begin{tikzpicture} [rotate=118]

  \draw \stA (0,1) -- (1,0) -- (-1,-1) -- cycle;
  \node at (0,0) {$\cdot$};

  \end{tikzpicture}

    $v_1=(1,0,1), v_2=(0,1,1), v_3=(-1,-1,1)$, $
    \det (v_1, v_2, v_3)  = 3$.
    This variety has a $\frac{1}{3}(1,1,1)$ quotient singularity, that is a quotient of an action by $\mathbb Z_3$.

    \item \textbf{Polygon P4:}  
      \begin{tikzpicture}[rotate=135]
 \node at (0,0) {$\cdot$};
  \draw \stA (0,1) -- (1,0) -- (-2,-1) -- (-1,0) -- cycle;

  \end{tikzpicture}
    
    \(v_1 = (-2, -1, 1)\), \(v_2 = (1, 0, 1)\),  \(v_3 = (0, 1, 1)\), $     \det (v_1, v_2, v_3) = 4$.
   This variety has a cyclic quotient singularity of an action by $\mathbb Z_4$, it is  \(\frac{1}{4}(1, 1, 2)\),  a common example in the study of terminal Gorenstein toric singularities.

    \item \textbf{Polygon P13:}
  \begin{tikzpicture}[rotate=90]
  \node at (0,0) {$\cdot$};
  \draw \stA (1,-1) -- (-1,-1) -- (-1,3) -- cycle;
  \end{tikzpicture}
  
     $v_1=(-1,-1,1), v_2=(1,-1,1), v_3=(-1,3,1)$, $     \det (v_1, v_2, v_3) = 8$.
    This variety has a quotient singularity by a group of order 8, in fact, by $\mathbb Z_2\times \mathbb Z_4$.

    \item \textbf{Polygon P16:} 
      \begin{tikzpicture}
  \draw \stA (2,-1) -- (-1,-1) -- (-1,2) -- cycle;
  \node at (0,0) {$\cdot$};
  \end{tikzpicture}
    
    $v_1=(-1,-1,1), v_2=(2,-1,1), v_3=(-1,2,1)$, $     \det (v_1, v_2, v_3) = 9$.
    This variety has a  quotient singularity by a group of order 9, in fact $\mathbb Z_3\times \mathbb Z_3$.

\end{itemize}

\newpage 

\section{Universal Quilts, Stacks, and Moduli}
\vspace{1cm}

Let $X$ be a singular toric Calabi--Yau threefold  (thus $X$ is necessarily noncompact) and let $P\ce P(X)$ be its corresponding polygon. 
Hence, if $(x_i,y_i)$ are the vertices of $P$, then $X$ is generated by the fan consisting of the set of  vectors that start at the origin $(0,0,0)$
and have endpoints  $(x_i,y_i,1)$.

\begin{definition}\label{primi}
A triangulation $T(P)$ of $P$ is called {\bf primitive} if it consists of a union of triangles of area $1/2$.
\end{definition}

\begin{definition} Let $\{T_\delta(P), \delta \in \Delta\}$ be the set of all triangulations of $P=P(X)$.  
A triangulated connected planar polygon ${\mathcal Q}$  is called  a {\bf versal quilt} for $X$ 
if it can be obtained as a union $${\mathcal Q}\ce{\mathcal Q}(X)= \bigcup_{\delta\in \Delta}  T_\delta(P) $$
where each triangulated polygon is attached to another (or others) by 
identifying boundary edges with equal length. 
\end{definition}

We observe that any primitive triangulation of $X$ can be obtained from  a quilt $Q(X)$ via pullback
by a (possibly non-unique) embedding. 

We will often refer to a versal quilt as just a {\bf quilt}, while reserving the more 
fundamental concept of universality for those quilts  satisfying the additional requirement 
that each individual primitive triangulation may be obtained from it via pullback in a unique way. 

\begin{definition}\label{qstack} A versal quilt ${\mathcal Q}$ for $X$ is called a  {\bf universal quilt} for $X$, denoted  $\mathcal  U(X)$
(or just $\mathcal U$ when clear from context), 
if it is universal for primitive triangulations of $P(X)$ with inclusions, that is, if each primitive triangulation of $X$ may be 
obtained from $\mathcal U$ via pullback by a unique map. 
\end{definition}

\begin{definition}
A {\bf quilt stack} for $X$ is a quotient  ${\mathcal Q}(X)/\sim$
where ${\mathcal Q}(X)$ is a versal quilt for $X$ and $\sim$  denotes equivalence 
via linear transformations of the planar polytopes by $SL^{\pm}(2,\mathbb Z)$  matrices. \footnote{We are 
grateful to Valerio Capellini for giving us the idea of considering equivalence classes modulo rotations and reflections.}
\end{definition}

We will represent quilt stacks using 3-dimensional 
figures,   where 
each individual primitive triangulation is represented by a height $1$ figure, 
and figures representing isomorphic  triangulations  are stacked as pilled up on a single column. 
Hence, the stalk of a 
primitive triangulation consists of all the elements on its  equivalence class modulo transformations 
by integer matrices with determinant $\pm1$, 
with the dimension of the stalk being its height.

Universal quilts can often be obtained from quilts by adding a finite number (possibly zero) of base points and marked points. 
This happens because adding base points
rigidifies the problem, imposing the condition that a base point of the primitive triangulation must be taken to a marked point
in the quilt. Next, we give examples of quilts, universal quilts and quilt stacks using our favourite polygons $Pi$  (notation \ref{polys}).

\clearpage
\subsection{Examples of Quilts and Stacks}\label{pictures}

\begin{example} The figure below represents a  versal quilt for polytope P16. Hence, the figure contains
\textbf{all 79 crepant resolutions of the Calabi--Yau threefold  $\mathbb C^3 / \mathbb Z_3 \times \mathbb Z_3.$}
Taking isomorphisms classes of its elements gives rise to the stack in figure \ref{stack1}.
\end{example}
\begin{center}
  \includegraphics[
  width=.94\textwidth,
  angle=180,
  trim={4mm 8mm 4mm 8mm},
  clip
  ]{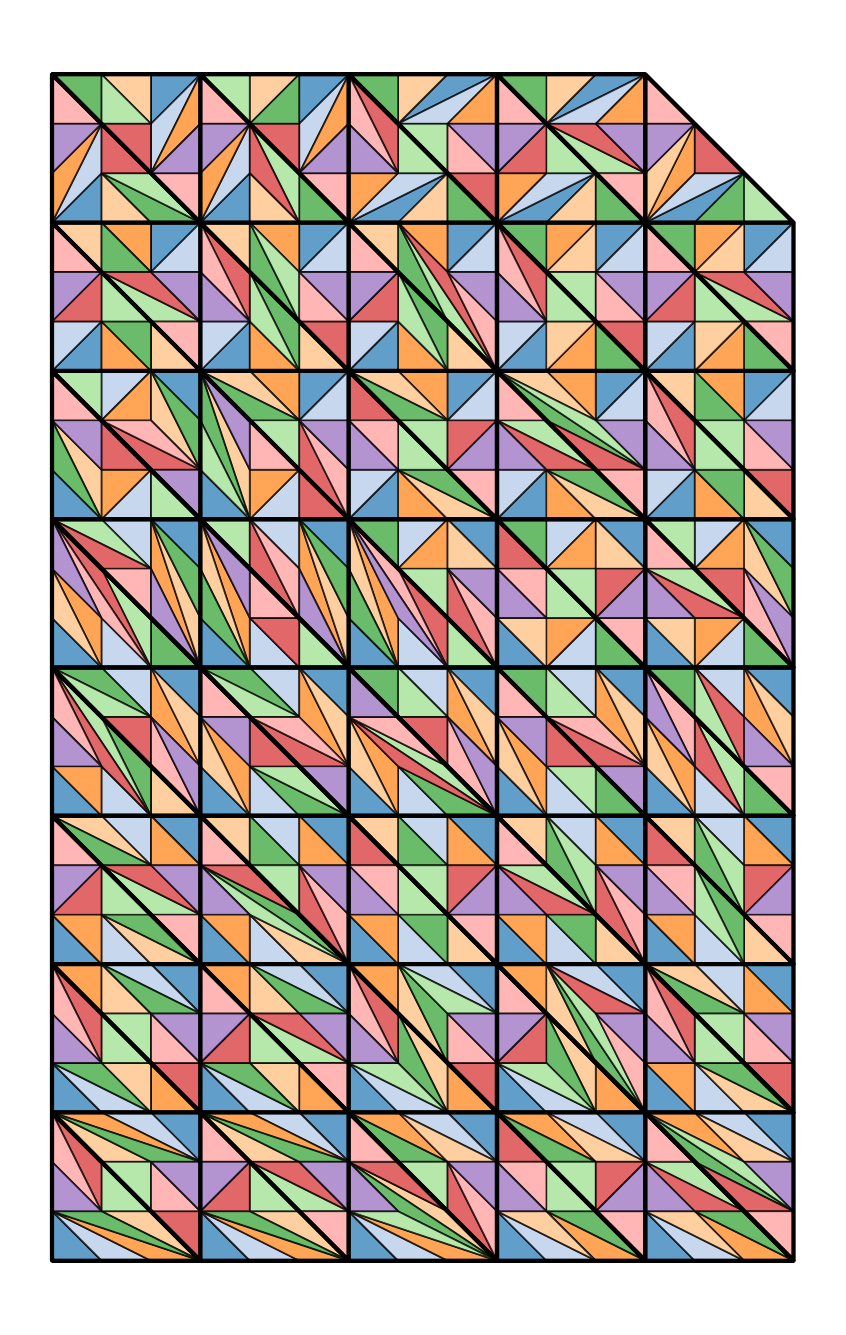}
\end{center}

\clearpage

\begin{example} The figure below represents a  universal quilt for polytope P16. Hence, the figure contains
\textbf{all 79 crepant resolutions of the Calabi--Yau threefold  $\mathbb C^3 / \mathbb Z_3 \times \mathbb Z_3$}
with basepoints, each embedded in the quilt  uniquely.
\end{example}
\begin{center}
  \includegraphics[
  width=.94\textwidth,
  trim={4mm 8mm 4mm 8mm},
  clip
  ]{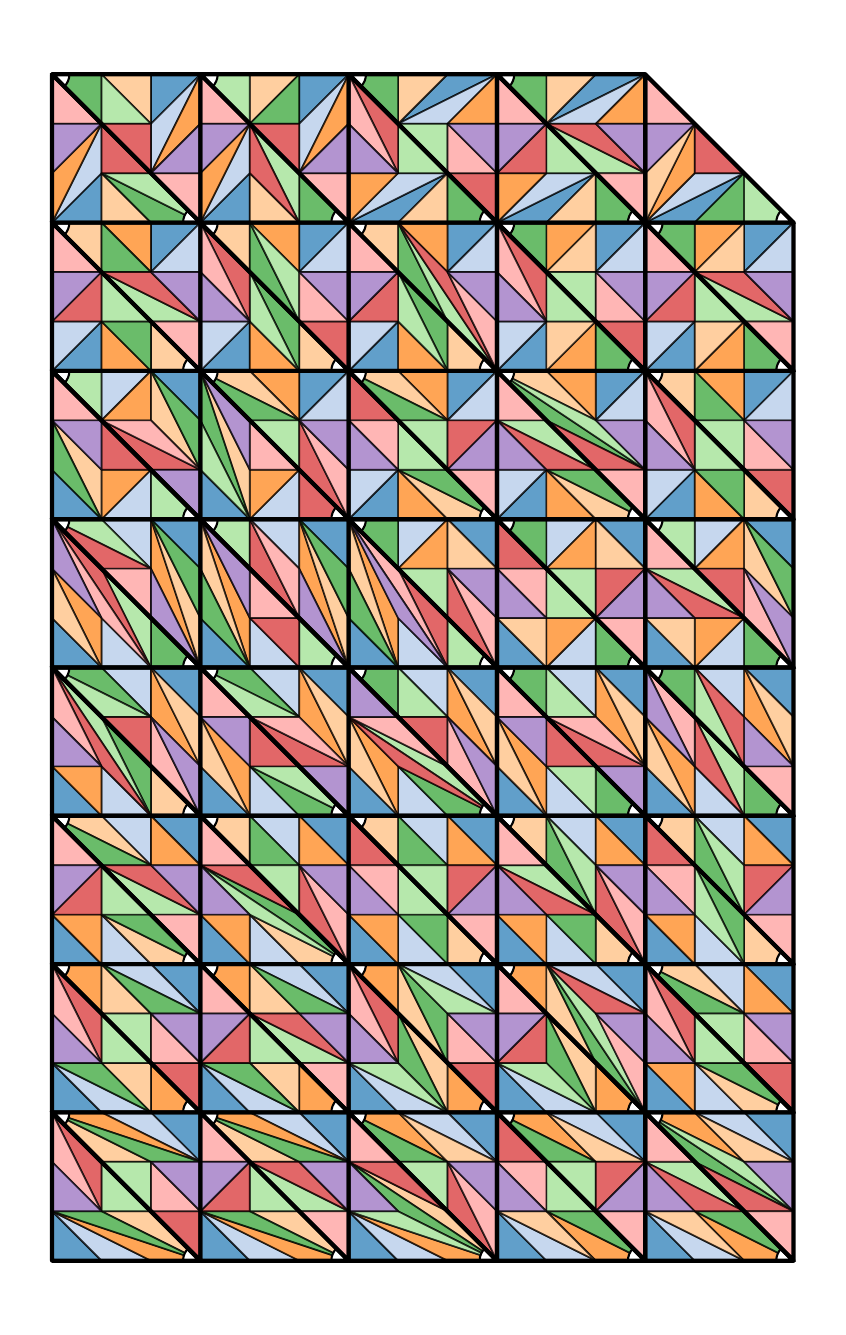}
\end{center}

\clearpage

\begin{example}\label{chess}
The figure below exemplifies a universal quilt for polytope P15. Thus, the quilt below
contains 
\textbf{all 64 crepant resolutions of the Calabi--Yau singularity determined by polygon P15}  
with basepoints. 
\end{example}
\begin{center}
  \includegraphics[
  width=.99\textwidth,
  trim={4mm 8mm 4mm 8mm},
  clip
  ]{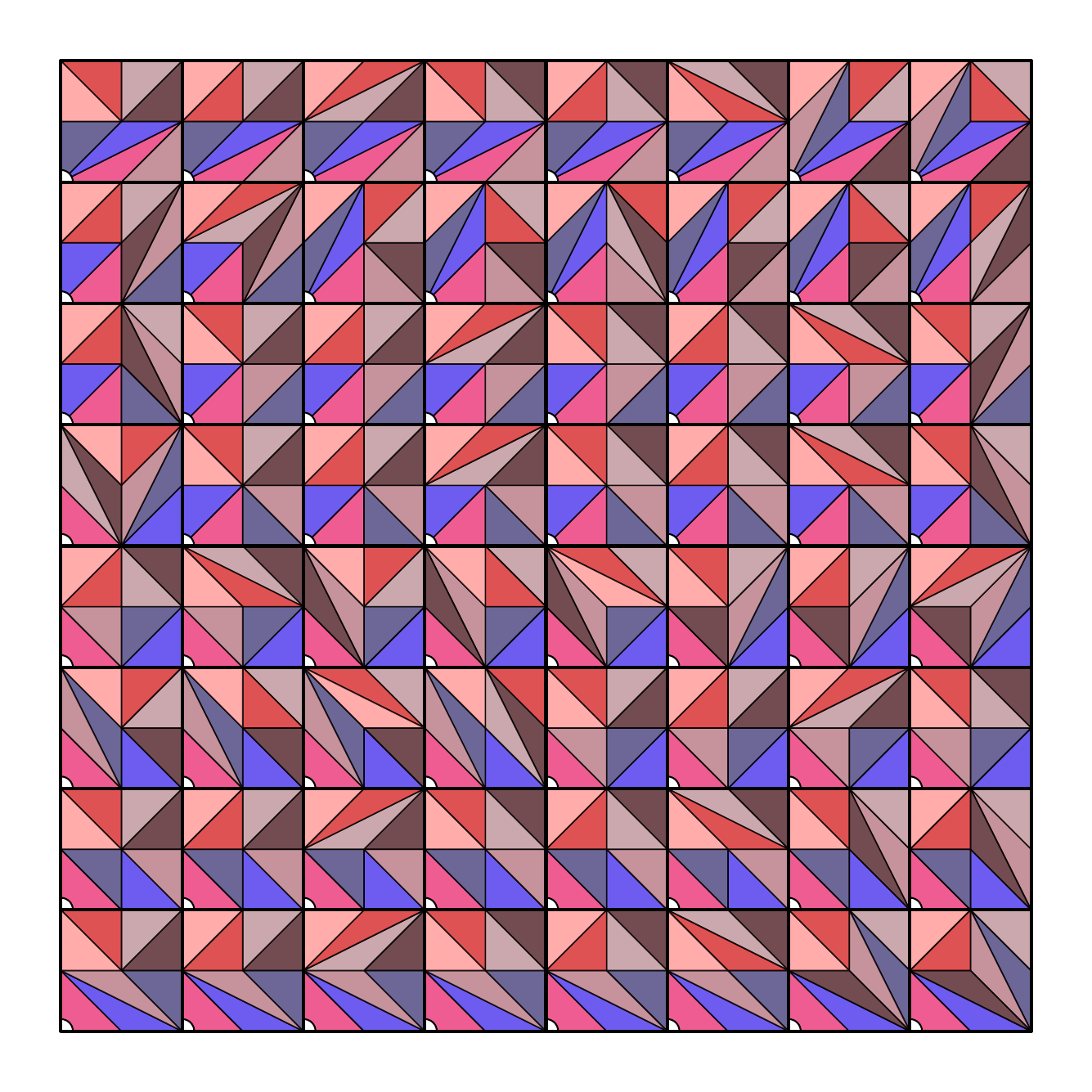}
\end{center}

We observe that this quilt was depicted in an 8x8 format, which happens to 
be the number of squares in a chess board (to which purpose a presentation 
without painting  half of the squares would be preferable).
Such a picture would still depict a universal quilt, because the colours here  
are given just for artistic illustration, i.e. just for fun.

\begin{minipage}{0.45\textwidth}
	\begin{example}\label{scarf}
	An alternative quilt for polygon P15. It contains the same triangulations as Example \ref{chess}, arranged and coloured differently.
\end{example}
\end{minipage}
\hfill
\begin{minipage}{0.45\textwidth}

\begin{center}
  \includegraphics[
	trim={0mm 0mm 0mm 0mm},
	scale=0.6,
	clip, 
	angle=90,
  ]{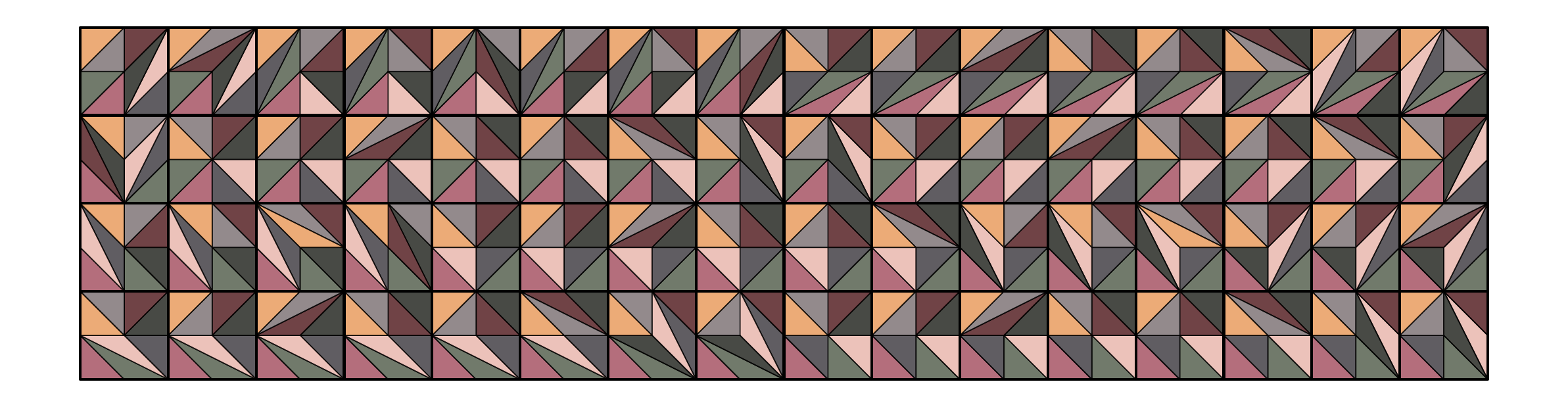}
\end{center}

\end{minipage}

\pagebreak

\begin{example} The figure below is a universal quit for polygon P7. Therefore, it has
\textbf{all 24 crepant resolutions of the Calabi--Yau threefold  $\mathbb C^3 / \mathbb Z_2 \times \mathbb Z_4.$}
\end{example}
\begin{center}
  \includegraphics[
  height=.88\textwidth,
  angle=90,
  trim={4mm 8mm 4mm 8mm},
  clip
  ]{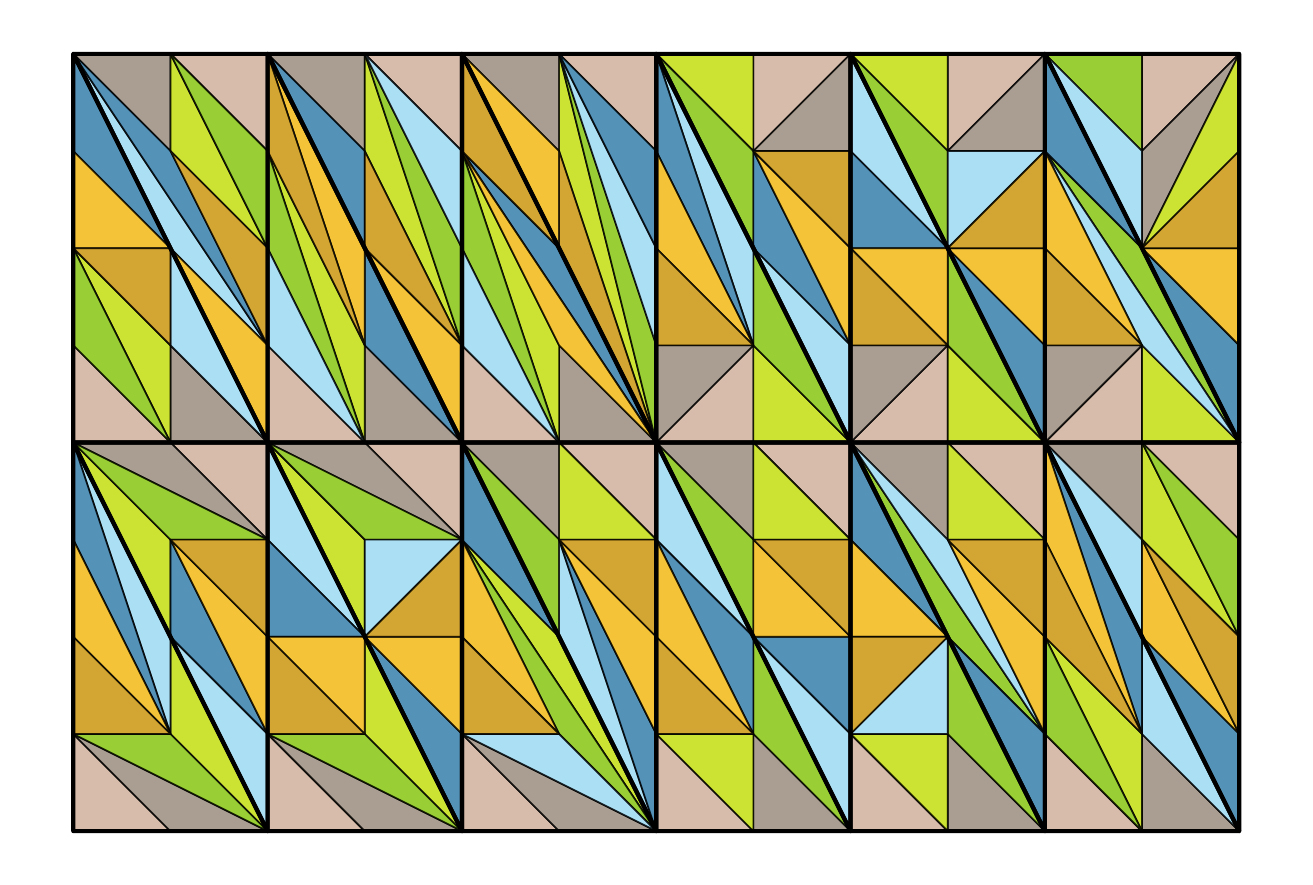}
\end{center}
\bigskip

We observe here no base points were required to obtain universality, and 
that this quilt also determines a fixed quilt stack.
That is, for polygon P13, the quilt and the stack are essentially the same.
This happens 
because in this case the  triangulations are pairwise non-isomorphic.
Hence, equivalence classes modulo $SL^{\pm}(2)$ transformations consist
of sets with a single element. 
Therefore, the stack would have maximum height 1.

\pagebreak

\begin{example}\label{ali}
\textbf{First  quilt for polygon P7}
\begin{center}
  \includegraphics[
  width=.60\textwidth,
  ]{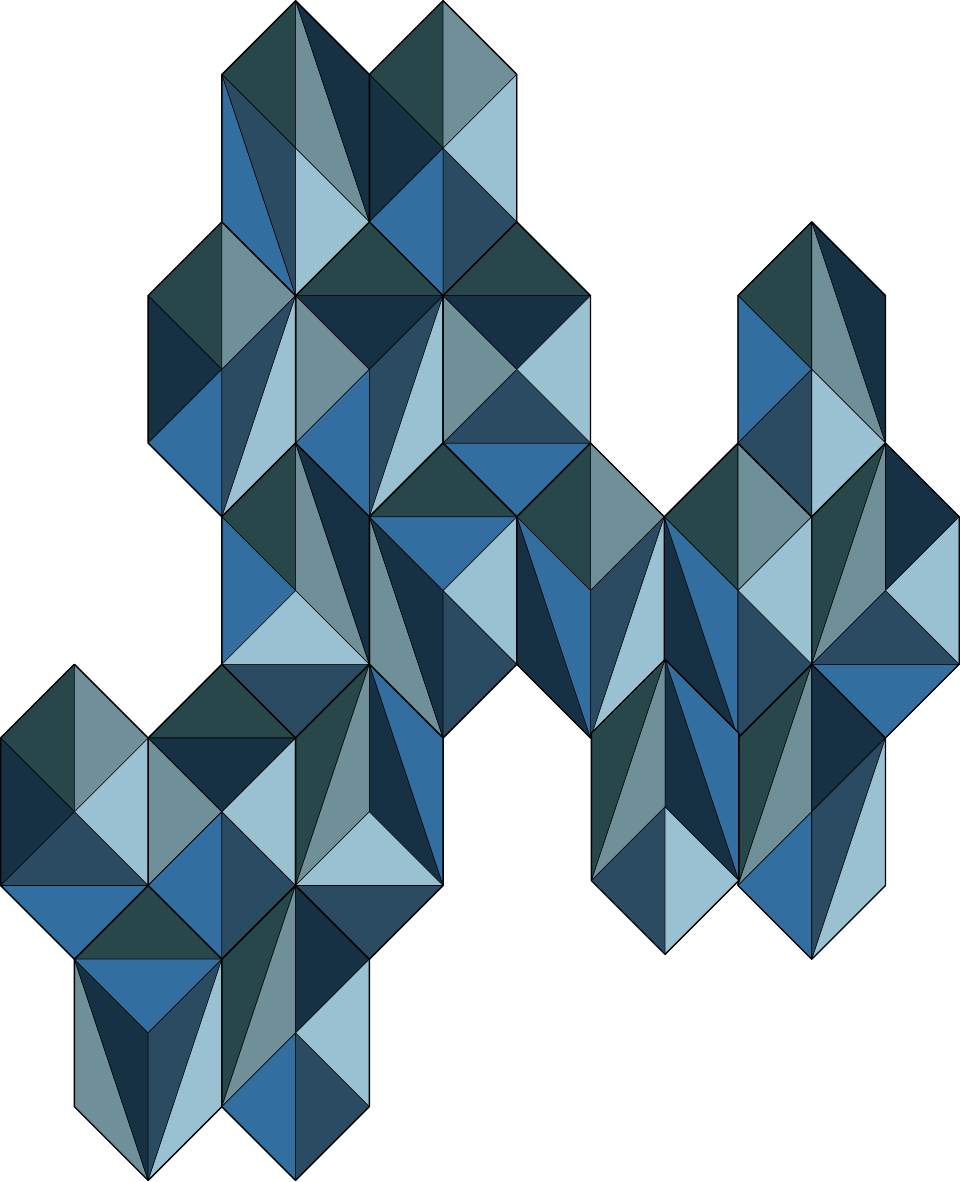}
\end{center}

\textbf{Second  quilt for polygon P7}
\begin{center}
  \includegraphics[
  width=.74\textwidth,
  angle=90,
  ]{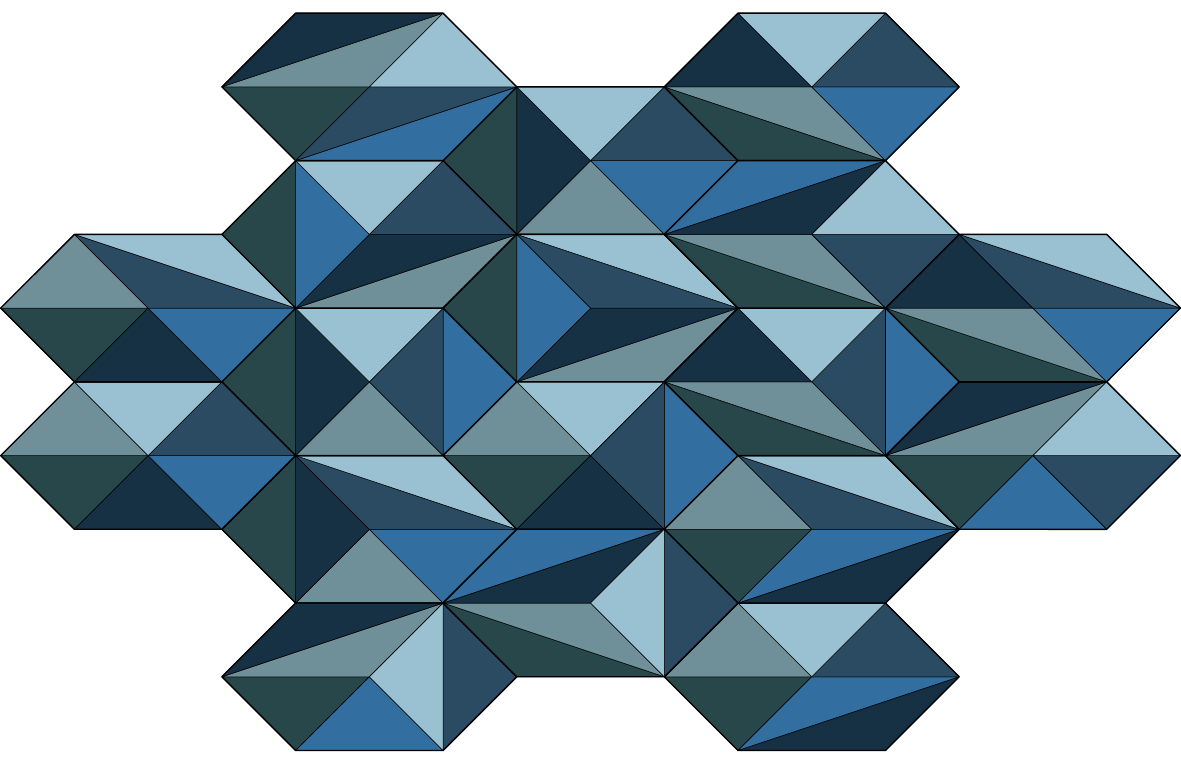}
\end{center}

The goals of this page are:
\begin{itemize}
\item  to exemplify the non-uniqueness of quilts, and
\item to emphasise the non-convexity of quilts for the hexagon (polygon P7).
\end{itemize}
\end{example}

\pagebreak

\begin{example} \label{stack1} We describe the moduli stack for polygon P15, with 
 the height of each column displaying the number of elements of the 
  corresponding $SL(3,\mathbb Z)$ isomorphism class. In this case the maximum height is 8.

 \begin{figure}[h]
 \begin{center}\includegraphics[height=10cm]{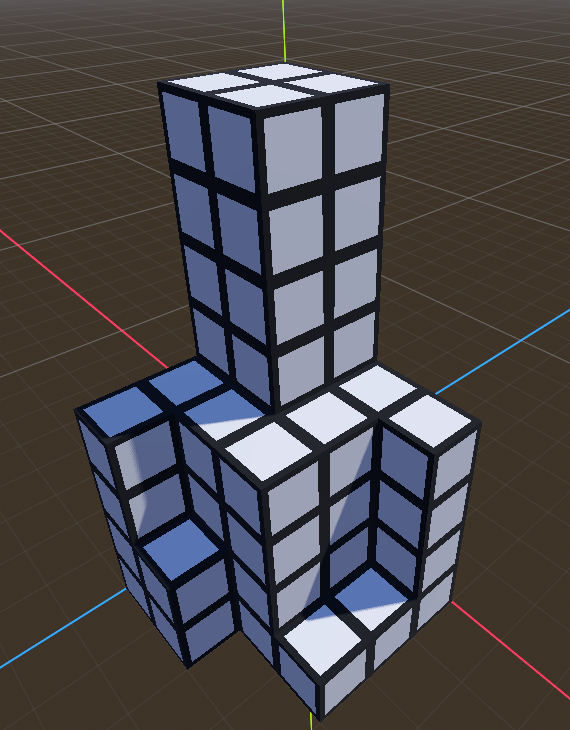}\end{center}
 \caption*{Moduli stack for polygon P15}
 \end{figure}
 
Here are 2 examples of triangulations whose equivalence class contains

a) just 1 element: 
 \vspace{3mm}
 
 \quad 
\begin{tikzpicture}[scale=\newscale]
\draw   (2,2) -- (2,0) -- (0,0)--(0,2)--(2,2);
\draw[green]  (2,1) -- (1,2);
\draw[green]  (0,1) -- (1,0);
\draw[green]  (2,1) -- (1,0);
\draw[green]  (0,1) -- (1,2);
\draw[red]  (2,1) -- (0,1);
\draw[red]  (1,2) -- (1,0);
\end{tikzpicture}

b) 8 elements:

\quad  
\begin{tikzpicture}[
scale=\newscale,
rotate=0
]
\begin{scope}[shift={(-1,-1)}]
\draw   (0,0) -- (2,0) -- (2,2) -- (0,2) -- cycle;
\draw[green]  (2,1) -- (1,2);
\draw[red]  (0,1) -- (1,0);
\draw[red]  (1,1) -- (2,0);
\draw[green]  (0,1) -- (2,0);
\draw[green]  (1,1) -- (0,2);
\draw[green]  (2,1) -- (1,1);
\draw[red]  (0,1) -- (1,1);
\draw[green]  (1,2) -- (1,1);
\end{scope}
\end{tikzpicture} 
\hfill
\begin{tikzpicture}[
scale=\newscale,
rotate=90
]
\begin{scope}[shift={(-1,-1)}]
\draw   (0,0) -- (2,0) -- (2,2) -- (0,2) -- cycle;
\draw[green]  (2,1) -- (1,2);
\draw[red]  (0,1) -- (1,0);
\draw[red]  (1,1) -- (2,0);
\draw[green]  (0,1) -- (2,0);
\draw[green]  (1,1) -- (0,2);
\draw[green]  (2,1) -- (1,1);
\draw[red]  (0,1) -- (1,1);
\draw[green]  (1,2) -- (1,1);
\end{scope}
\end{tikzpicture} 
\hfill
\begin{tikzpicture}[
scale=\newscale,
rotate=180
]
\begin{scope}[shift={(-1,-1)}]
\draw   (0,0) -- (2,0) -- (2,2) -- (0,2) -- cycle;
\draw[green]  (2,1) -- (1,2);
\draw[red]  (0,1) -- (1,0);
\draw[red]  (1,1) -- (2,0);
\draw[green]  (0,1) -- (2,0);
\draw[green]  (1,1) -- (0,2);
\draw[green]  (2,1) -- (1,1);
\draw[red]  (0,1) -- (1,1);
\draw[green]  (1,2) -- (1,1);
\end{scope}
\end{tikzpicture} 
\hfill
\begin{tikzpicture}[
scale=\newscale,
rotate=270
]
\begin{scope}[shift={(-1,-1)}]
\draw   (0,0) -- (2,0) -- (2,2) -- (0,2) -- cycle;
\draw[green]  (2,1) -- (1,2);
\draw[red]  (0,1) -- (1,0);
\draw[red]  (1,1) -- (2,0);
\draw[green]  (0,1) -- (2,0);
\draw[green]  (1,1) -- (0,2);
\draw[green]  (2,1) -- (1,1);
\draw[red]  (0,1) -- (1,1);
\draw[green]  (1,2) -- (1,1);
\end{scope}
\end{tikzpicture} 
\hfill
\begin{tikzpicture}[
scale=\newscale,
xscale=-1,
rotate=0
]
\begin{scope}[shift={(-1,-1)}]
\draw   (0,0) -- (2,0) -- (2,2) -- (0,2) -- cycle;
\draw[green]  (2,1) -- (1,2);
\draw[red]  (0,1) -- (1,0);
\draw[red]  (1,1) -- (2,0);
\draw[green]  (0,1) -- (2,0);
\draw[green]  (1,1) -- (0,2);
\draw[green]  (2,1) -- (1,1);
\draw[red]  (0,1) -- (1,1);
\draw[green]  (1,2) -- (1,1);
\end{scope}
\end{tikzpicture}
\hfill
\begin{tikzpicture}[
scale=\newscale,
xscale=-1,
rotate=90
]
\begin{scope}[shift={(-1,-1)}]
\draw   (0,0) -- (2,0) -- (2,2) -- (0,2) -- cycle;
\draw[green]  (2,1) -- (1,2);
\draw[red]  (0,1) -- (1,0);
\draw[red]  (1,1) -- (2,0);
\draw[green]  (0,1) -- (2,0);
\draw[green]  (1,1) -- (0,2);
\draw[green]  (2,1) -- (1,1);
\draw[red]  (0,1) -- (1,1);
\draw[green]  (1,2) -- (1,1);
\end{scope}
\end{tikzpicture}
\hfill
\begin{tikzpicture}[
scale=\newscale,
xscale=-1,
rotate=180
]
\begin{scope}[shift={(-1,-1)}]
\draw   (0,0) -- (2,0) -- (2,2) -- (0,2) -- cycle;
\draw[green]  (2,1) -- (1,2);
\draw[red]  (0,1) -- (1,0);
\draw[red]  (1,1) -- (2,0);
\draw[green]  (0,1) -- (2,0);
\draw[green]  (1,1) -- (0,2);
\draw[green]  (2,1) -- (1,1);
\draw[red]  (0,1) -- (1,1);
\draw[green]  (1,2) -- (1,1);
\end{scope}
\end{tikzpicture}
\hfill
\begin{tikzpicture}[
scale=\newscale,
xscale=-1,
rotate=270
]
\begin{scope}[shift={(-1,-1)}]
\draw   (0,0) -- (2,0) -- (2,2) -- (0,2) -- cycle;
\draw[green]  (2,1) -- (1,2);
\draw[red]  (0,1) -- (1,0);
\draw[red]  (1,1) -- (2,0);
\draw[green]  (0,1) -- (2,0);
\draw[green]  (1,1) -- (0,2);
\draw[green]  (2,1) -- (1,1);
\draw[red]  (0,1) -- (1,1);
\draw[green]  (1,2) -- (1,1);
\end{scope}
\end{tikzpicture}
\vspace{2mm}

So, the conclusion is that  triangulation b) above is 8 times more likely to occur 
than triangulation a).

Observe that the base of this moduli stack consists of 14 boxes, which is the number of 
moduli as described in example \ref{trigmod}. Our stack figures are very similar to the figures appearing in \cite{ORV}.
\end{example}
 \clearpage
 
\begin{example} \label{stack2} We describe the moduli stack for polygon 7, with 
 the height of each column displaying the number of elements of the 
  corresponding $SL(3,\mathbb Z)$ isomorphism class. In this case the maximum height is 6.

 \begin{figure}[h]
 \begin{center}\includegraphics[height=6cm]{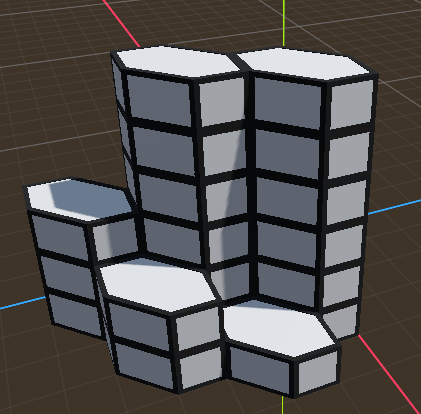}\end{center}
 \caption*{Moduli stack for polygon P7}
 \end{figure}

Here are examples of triangulations whose equivalence class contains

a) only 1 element:
 \vspace{3mm}
 
  \begin{tikzpicture}[scale=\newscale]
  \draw[green] (0,0) -- (0,1);
  \draw[green]  (0,0) -- (1,1);
  \draw[green]  (0,0) -- (1,0);
  \draw[green]  (0,0) -- (0,-1);
  \draw[green]  (0,0) -- (-1,-1);
  \draw[green]  (0,0) -- (-1,0);

  \draw  (0,1) -- (1,1) -- (1,0) -- (0,-1) -- (-1,-1) -- (-1,0) -- cycle;

  \end{tikzpicture}
 \vspace{3mm}
 
b) 6 elements:

 \vspace{3mm}

\begin{tikzpicture}[
scale=\newscale,
rotate=0
]
\begin{scope}[shift={(0,0)}]
\draw[red] (0,0) -- (0,1);
\draw[green]  (0,0) -- (1,1);
\draw[red]  (0,0) -- (1,0);
\draw[green]  (-1,-1) -- (1,0);
\draw[cyan]  (0,0) -- (-1,-1);
\draw[green]  (-1,-1) -- (0,1);

\draw  (0,1) -- (1,1) -- (1,0) -- (0,-1) -- (-1,-1) -- (-1,0) -- cycle;
\end{scope}
\end{tikzpicture}
\hfill
\begin{tikzpicture}[
scale=\newscale,
rotate=180
]
\begin{scope}[shift={(0,0)}]
\draw[red] (0,0) -- (0,1);
\draw[green]  (0,0) -- (1,1);
\draw[red]  (0,0) -- (1,0);
\draw[green]  (-1,-1) -- (1,0);
\draw[cyan]  (0,0) -- (-1,-1);
\draw[green]  (-1,-1) -- (0,1);

\draw  (0,1) -- (1,1) -- (1,0) -- (0,-1) -- (-1,-1) -- (-1,0) -- cycle;
\end{scope}
\end{tikzpicture}
\hfill
\begin{tikzpicture}[
scale=\newscale,
xscale=1,
rotate=0
]
\begin{scope}[shift={(0,0)}]
\draw[red] (0,0) -- (1,0);
\draw[red] (0,0) -- (-1,-1);
\draw[green] (0,0) -- (0,-1);
\draw[green] (1,0) -- (0,1);
\draw[green] (-1,-1) -- (0,1);
\draw[blue] (0,0) -- (0,1);
\draw  (0,1) -- (1,1) -- (1,0) -- (0,-1) -- (-1,-1) -- (-1,0) -- cycle;
\end{scope}
\end{tikzpicture}
\hfill
\begin{tikzpicture}[
scale=\newscale,
xscale=-1,
rotate=90
]
\begin{scope}[shift={(0,0)}]
\draw[red] (0,0) -- (1,0);
\draw[red] (0,0) -- (-1,-1);
\draw[green] (0,0) -- (0,-1);
\draw[green] (1,0) -- (0,1);
\draw[green] (-1,-1) -- (0,1);
\draw[blue] (0,0) -- (0,1);
\draw  (0,1) -- (1,1) -- (1,0) -- (0,-1) -- (-1,-1) -- (-1,0) -- cycle;
\end{scope}
\end{tikzpicture}
\hfill
\begin{tikzpicture}[
scale=\newscale,
xscale=1,
rotate=0
]
\begin{scope}[shift={(0,0)}]
\draw[red] (0,0) -- (1,1);
\draw[red] (0,0) -- (0,-1);
\draw[green] (-1,0) -- (1,1);
\draw[green] (-1,0) -- (0,-1);
\draw[green] (0,0) -- (1,0);
\draw[blue] (0,0) -- (-1,0);
\draw  (0,1) -- (1,1) -- (1,0) -- (0,-1) -- (-1,-1) -- (-1,0) -- cycle;
\end{scope}
\end{tikzpicture}
\hfill
\begin{tikzpicture}[
scale=\newscale,
xscale=-1,
rotate=90
]
\begin{scope}[shift={(0,0)}]
\draw[red] (0,0) -- (1,1);
\draw[red] (0,0) -- (0,-1);
\draw[green] (-1,0) -- (1,1);
\draw[green] (-1,0) -- (0,-1);
\draw[green] (0,0) -- (1,0);
\draw[blue] (0,0) -- (-1,0);
\draw  (0,1) -- (1,1) -- (1,0) -- (0,-1) -- (-1,-1) -- (-1,0) -- cycle;
\end{scope}
\end{tikzpicture}
 \vspace{3mm}  
  
So, the conclusion is that  triangulation b) above is 6 times more likely to occur 
than triangulation a).

Observe that the base of the moduli stack consists of 6 boxes, which is the number of 
moduli as described in example \ref{trigmod}.

\end{example}

\clearpage
\section{Theorems}

If a variety $X$ admits a convex universal quilt $\mathcal U(X)$, then such a quilt  itself represents a 
smooth Calabi--Yau threefold, given that it comes by construction endowed with a primitive triangulation.

\begin{theorem}\label{tri}
Let $X$ be a toric Calabi--Yau threefold singularity whose planar polytope 
is a convex triangle. Then it admits a smooth toric Calabi--Yau threefold as a universal quilt $\mathcal U(X)$. 
\end{theorem}

\begin{theorem}\label{par}
Let $X$ be a toric Calabi--Yau threefold singularity whose planar polytope 
is a parallelogram. Then it admits a smooth toric Calabi--Yau threefold as a universal quilt $\mathcal U(X)$. 
\end{theorem}

\begin{theorem}\label{tah}
Let $X$ be a toric Calabi--Yau threefold singularity whose planar polytope 
is a hexagon. Then the quilts of $X$ do not represent toric varieties. 
\end{theorem}

\paragraph{\bf Acknowledgements}  
This work was supported by the European Research Council through the Horizon ERC Synergy Grant 101167526 (MALINCA).
E. Gasparim is a Senior Associate of the Abdus Salam International Centre for Theoretical Physics, Trieste, and is grateful for insightful conversations at ICTP; in particular with Micol Stock, Emanuel Carneiro, Gabriella Orlando, and Cumrun Vafa.
R. Talha gratefully acknowledges the support from the ICTP Postgraduate Diploma Programme during his studies.

\end{document}